\documentclass[11pt]{article}
\usepackage{graphicx,wrapfig,subfigure}
\usepackage{amsmath}
\usepackage{amsthm}
 \usepackage{amssymb}
\usepackage{setspace}
\usepackage{leftidx}
\usepackage{graphicx}
\usepackage{amsfonts}
\usepackage{color, epstopdf}
\usepackage{cite}
\usepackage{array,setspace,diagbox}
\usepackage{booktabs}
\usepackage [latin1]{inputenc}
\usepackage{multirow}
\usepackage{enumitem}
\usepackage{algpseudocode,algorithm,algorithmicx}%
\usepackage{appendix}%
\usepackage{upgreek}
\usepackage{enumerate}
\usepackage{enumitem}

\newcommand{\tabincell}[2]{\begin{tabular}{@{}#1@{}}#2\end{tabular}}

\newtheorem{theorem}{Theorem}[section]
\newtheorem{lemma}{Lemma}[section]

\newtheorem{corollary}{Corollary}[section]
\newtheorem{example}{Example}
\newtheorem{remark}{Remark}

\newcommand{\diff}{\triangledown_{\!\tau}}

\newcommand{\zd}{\,\mathrm{d}}
\newcommand{\abs}[1]{\left|#1\right|}
\newcommand{\absb}[1]{\big|#1\big|}

\newcommand{\bra}[1]{\left(#1\right)}
\newcommand{\brab}[1]{\big(#1\big)}
\newcommand{\braB}[1]{\Big(#1\Big)}
\newcommand{\brat}[1]{(#1)}
\newcommand{\kbra}[1]{\left[#1\right]}
\newcommand{\kbrab}[1]{\big[#1\big]}
\newcommand{\kbraB}[1]{\Big[#1\Big]}

\newcommand{\mynorm}[1]{\left\|#1\right\|}
\newcommand{\mynormb}[1]{\big\|#1\big\|}
\newcommand{\mynormB}[1]{\Big\|#1\Big\|}

\newcommand{\Lh}[2]{a_{#1}^{(\mathrm{h},#2)}} 
\newcommand{\La}[2]{a_{#1}^{(\mathrm{a},#2)}} 
\newcommand{\aLh}[3]{\mathsf{#1}_{#2}^{(\mathrm{h},#3)}} 
\newcommand{\aLa}[3]{\mathsf{#1}_{#2}^{(\mathrm{a},#3)}} 

\def\lan#1{\textcolor{blue}{#1}}

\begin{document}
\title{Energy stability of variable-step L1-type schemes
 for time-fractional Cahn-Hilliard model}
\author{
	Bingquan Ji\,\footnotemark[1] \footnotemark[2]
	\quad Xiaohan Zhu\,\footnotemark[1] \footnotemark[2]
	\quad Hong-lin Liao\,\footnotemark[3]}
\footnotetext[1]{$^{*}$ Department of Mathematics,
	Nanjing University of Aeronautics and Astronautics, Nanjing 211106, P. R. China.}
\footnotetext[2]{$^{\dag}$ These authors contributed equally to this work.}
\footnotetext[3]{$^{\ddag}$ Corresponding author. ORCID 0000-0003-0777-6832. Department of Mathematics,
	Nanjing University of Aeronautics and Astronautics,
	Nanjing 211106, China; Key Laboratory of Mathematical Modelling
	and High Performance Computing of Air Vehicles (NUAA), MIIT, Nanjing 211106, China.
	Hong-lin Liao (liaohl@nuaa.edu.cn and liaohl@csrc.ac.cn)
	is supported by a grant 12071216 from
	National Natural Science Foundation of China.}
\date{\today}
\maketitle
\normalsize

\begin{abstract}
The positive definiteness of discrete time-fractional derivatives is fundamental
to the numerical stability (in the energy sense) for time-fractional phase-field models.
A novel technique is proposed to estimate the minimum eigenvalue
of discrete convolution kernels generated by the nonuniform  L1, half-grid based L1 and time-averaged L1 formulas
of the fractional Caputo's derivative. The main discrete tools are
the discrete orthogonal convolution kernels and discrete complementary convolution kernels.
Certain variational energy dissipation laws at discrete levels of the variable-step
 L1-type methods are then established for time-fractional Cahn-Hilliard model.
They are shown to be asymptotically compatible,
in the fractional order limit $\alpha\rightarrow1$, with the associated energy dissipation law for the classical Cahn-Hilliard equation.
Numerical examples together with an adaptive time-stepping procedure
 are provided to demonstrate the effectiveness of the proposed methods.\\
  \noindent{\emph{Keywords}:}\;\; time-fractional Cahn-Hilliard model;
  variable-step L1-type  formulas; discrete convolution tools;
  positive definiteness; variational energy dissipation law\\
  \noindent{\bf AMS subject classiffications.}\;\; 35Q99, 65M06, 65M12, 74A50
\end{abstract}

\section{Introduction}
\setcounter{equation}{0}

In this paper, we study the energy stability of three nonuniform L1-type approximations
for the time-fractional Cahn-Hilliard (TFCH) model\cite{TangYuZhou:2018TimeFractionalPhase}
\begin{align}\label{cont:TFCH model}
\partial_t^\alpha \Phi=\kappa\Delta\mu\quad\text{with}\quad
\mu=\tfrac{\delta E}{\delta \Phi}=f(\Phi)-\epsilon^2\Delta\Phi,
\end{align}
where the Ginzburg-Landau energy functional
is given by \cite{CahnHilliard:1958Free},
\begin{align}\label{cont:classical free energy}
E[\Phi]=\int_{\Omega}\braB{\frac{\epsilon^2}{2}\abs{\nabla\Phi}^2
+F(\Phi)}\zd\mathbf{x}.
\end{align}
Here, the real valued function $\Phi$ represents the concentration
difference in a binary system,
spatial domain $\mathbf{x}\in\Omega\subseteq\mathbb{R}^2$,
$\epsilon>0$ is an interface width parameter, \lan{ $\kappa>0$ is the mobility coefficient}
and the double-well potential $F(\Phi)=\frac{1}{4}\bra{\Phi^2-1}^2$.
The notation $\partial_t^\alpha:={}_{0}^{C}\!D_{t}^{\alpha}$ represents
the Caputo's fractional  derivative of order $\alpha\in(0,1)$
with respect to $t$, defined by \cite{Podlubny:1999FractionalBook},
\begin{align}\label{def:Caputo derivative}
(\partial_{t}^{\alpha}v)(t)
:=(\mathcal{I}_{t}^{1-\alpha}v^\prime)(t)
\quad \text{where}\quad
(\mathcal{I}_{t}^{\beta}v)(t)
:=\int_{0}^{t}\omega_{\beta}(t-s)v(s)\zd{s},
\end{align}
and $\omega_{\beta}(t):=t^{\beta-1}/\Gamma(\beta)$ for $\beta>0$.

Throughout this paper, the periodic boundary conditions are adopted for simplicity.
\lan{
If the initial data is properly regular,
the global existence of solutions of the TFCH equation \eqref{cont:TFCH model} was established in 
\cite{Karaa:2021TFCHGlobalEnergy}.
Moreover, \cite[Theorem 3.3]{Karaa:2021TFCHGlobalEnergy} showed that the problem \eqref{cont:TFCH model} has a 
unique solution and $\mynormb{\partial_t\Phi}\le C_{\Phi} t^{\alpha\nu/4-1}$ for $0<t \le T$ if $\Phi_0\in \dot{H}^{\nu}(\Omega)\,(\nu\in[1,2])$.
It reveals that
the solution of the TFCH equation lacks the smoothness near the initial time while it would be smooth away from $t=0$,
also see \cite{Stynes:2017ErrorAnalysis,Kopteva:2019Error,LiaoYanZhang:2018Unconditional}
and the references therein.}
In addition,
the TFCH equation \eqref{cont:TFCH model}
conserves the initial volume \lan{ $\brab{\Phi(t),1}=\bra{\Phi(0),1}$ }for $t>0$
\cite[Theorem 2.2]{TangYuZhou:2018TimeFractionalPhase}.
Recently, Liao, Tang and Zhou \cite{LiaoTangZhou:2020EnergyMaximumBound} showed that
the time-fractional phase field models preserve
the following variational energy dissipation law,
\begin{align}\label{cont:TFCH variational energy}
\frac{\zd \mathcal{E}_\alpha}{\zd t}
+\frac{\kappa}{2}\omega_{\alpha}(t)\mynormb{\nabla\mu}^2\le0\quad
\text{with} \quad \mathcal{E}_\alpha(t):=E(t)
+\frac{\kappa}{2}\mathcal{I}_{t}^{\alpha}\mynormb{\nabla\mu}^2
\quad \text{for $t>0$},
\end{align}
where $(\cdot,\cdot)$ and $\mynorm{\cdot}$ denote the
$L^2(\Omega)$ inner product and the associated norm, respectively.

As pointed out in \cite{LiaoTangZhou:2020EnergyMaximumBound},
compared with the global energy dissipation property
in \cite{TangYuZhou:2018TimeFractionalPhase,ChenZhangZhao:2019TimeMBE},
the time-fractional energy decaying law
and the weighted energy dissipation law in \cite{QuanTangYang:22020HowToDefine,QuanTangYang:2020Numerical},
the new law \eqref{cont:TFCH variational energy} seems to be naturally
consistent with the standard energy dissipation law of the classical
Cahn-Hilliard (CH) model in the sense that
\begin{align*}
\frac{\zd \mathcal{E}_\alpha}{\zd t}
+\frac{\kappa}{2}\omega_{\alpha}(t)\mynormb{\nabla\mu}^2\le0\quad
\rightarrow\quad\frac{\zd E}{\zd t}+\kappa\mynormb{\nabla\mu}^2\le0
\qquad \text{as $\alpha\rightarrow 1$}.
\end{align*}

Our aim of this paper is to develop numerical methods
that preserve the variational energy dissipation law
\eqref{cont:TFCH variational energy} at discrete  time levels.
For a given time $T>0$,
consider a nonuniform time levels $0=t_{0}<t_{1}<\cdots<t_{k-1}<t_{k}<\cdots<t_{N}=T$
with the time-step sizes $\tau_{k}:=t_{k}-t_{k-1}$ for $1\le k\le N$.
The maximum time-step size is denoted by $\tau:=\max_{1\le k\le N}\tau_k$
and the local time-step ratio $r_k:=\tau_k/\tau_{k-1}$ for $2\le k\le N$.

Given a grid function $\{v^{k}\}_{k=0}^N$, define the difference
$\diff v^k:=v^k-v^{k-1}$ and
$\partial_\tau v^k:=\diff v^k/\tau_k$ for $k\geq{1}$.
Let $\Pi_{1,k}v$ denote the linear interpolant of a function $v$
with respect to the nodes $t_{k-1}$ and $t_k$, such that
$\bra{\Pi_{1,k}v}^\prime(t)=\diff v^k/\tau_k$ for $t\in (t_{k-1},t_k)$.
We will investigate three L1-type formulas on nonuniform meshes.
The first one is the standard L1 approximation
\cite{LiaoYanZhang:2018Unconditional},
 \begin{align}\label{scheme:L1 formula}
(\partial_{\tau}^{\alpha}v)^n
:=\int_0^{t_n}
\omega_{1-\alpha}(t_n-s)(\Pi_{1}v)^\prime(s)\zd{s}
\triangleq\sum_{k=1}^{n}a_{n-k}^{(n)}\diff v^{k},
\end{align}
where the associated discrete L1 kernels $a_{n-k}^{(n)}$ are defined by
\lan{
\begin{align}\label{def:L1 formula kernel}
a_{n-k}^{(n)}&:=\frac{1}{\tau_{k}}\int_{t_{k-1}}^{t_{k}}\omega_{1-\alpha}(t_n-s)\zd{s} \nonumber \\
&=\frac{1}{\tau_{k}}\kbra{\omega_{2-\alpha}(t_n-t_{k-1})-\omega_{2-\alpha}(t_n-t_{k})}
\quad\text{for $1\leq{k}\leq{n}$}.
\end{align}
}
The second formula, named L1$_{\mathrm{h}}$,  is defined at the half-grid point $t_{n-\frac12}$,
\begin{align}\label{scheme:L1h formula}
(\partial_{\mathrm{h}\tau}^{\alpha}v)^{n-\frac12}
:=\int_0^{t_{n-\frac12}}
\omega_{1-\alpha}(t_{n-\frac12}-s)(\Pi_{1}v)'(s)\zd{s}
\triangleq\sum_{k=1}^{n}\Lh{n-k}{n}\diff v^{k},
\end{align}
where the corresponding discrete L1$_\mathrm{h}$ kernels $\Lh{n-k}{n}$ are given by
\begin{align}\label{def:L1h formula kernel}
\Lh{n-k}{n}:=\frac{1}{\tau_{k}}\int_{t_{k-1}}^{\min\{t_k,t_{n-\frac12}\}}
\omega_{1-\alpha}(t_{n-\frac12}-s)\zd{s}
\quad\text{for $1\le k\le n$}.
\end{align}
The third one, called  L1$_\mathrm{a}$, is  an averaged version of L1 formula \eqref{scheme:L1 formula}
at  $t_{n-\frac12}$, that is,
\begin{align}\label{scheme:L1a formula}
(\partial_{\mathrm{a}\tau}^{\alpha}v)^{n-\frac12}
:=\frac12\kbra{(\partial_{\tau}^{\alpha}v)^n+(\partial_{\tau}^{\alpha}v)^{n-1}}
\triangleq\sum_{k=1}^{n}\La{n-k}{n}\diff v^{k},
\end{align}
where the corresponding discrete L1$_\mathrm{a}$ kernels $\La{n-k}{n}$
are defined by
\begin{align}\label{def:L1a formula kernel}
a_{0}^{(\mathrm{a},n)}:=\frac12a_{0}^{(n)}\quad\text{and}\quad
a_{n-k}^{(\mathrm{a},n)}
:=\frac12\brab{a_{n-k}^{(n)}+a_{n-1-k}^{(n-1)}}
\quad\text{for $1\le k\le n-1$}.
\end{align}

By means of the above three L1-type formulas
\eqref{def:L1h formula kernel}-\eqref{def:L1a formula kernel},
we consider the following semi-discrete time-stepping methods for the TFCH model:
\begin{itemize}[itemindent=0.5em]
\item The backward Euler-type L1 scheme
\begin{align}\label{scheme:fully implicit L1 TFCH}
\bra{\partial_\tau^\alpha \phi}^n=\kappa\Delta \mu^n\quad\text{with}\quad
\mu^n=f(\phi^n)-\epsilon^2\Delta \phi^n, \quad  n\ge 1;
\end{align}

\item The Crank-Nicolson-type L1$_\mathrm{h}$  scheme
\begin{align}\label{scheme:fully implicit L1h TFCH}
\bra{\partial_{\mathrm{h}\tau}^\alpha \phi}^{n-\frac12}
=\kappa\Delta \mu^{n-\frac12}\quad\text{with}\quad
\mu^{n-\frac12}=f\brat{\phi}^{n-\frac12}-\epsilon^2\Delta \phi^{n-\frac12},
\quad  n\ge 1;
\end{align}

\item The Crank-Nicolson-type L1$_\mathrm{a}$ scheme
\begin{align}\label{scheme:fully implicit L1a TFCH}
\bra{\partial_{\mathrm{a}\tau}^\alpha \phi}^{n-\frac12}
=\kappa\Delta \mu^{n-\frac12}\quad\text{with}\quad
\mu^{n-\frac12}=f\brat{\phi}^{n-\frac12}-\epsilon^2\Delta \phi^{n-\frac12},
\quad  n\ge 1,
\end{align}
\end{itemize}
where the averaged difference operator $\phi^{n-\frac12}:=(\phi^n+\phi^{n-1})/2$
and $f\brat{\phi}^{n-\frac12}$ is a second-order
approximation \cite[Appendix A]{LiaoTangZhou:2020EnergyMaximumBound}
 of the nonlinear term $f(\phi)$
\[
f\brat{\phi}^{n-\frac12}:=\frac13\bra{\phi^n}^3+\frac12\phi^n\bra{\phi^{n-1}}^2
+\frac16\bra{\phi^{n-1}}^3-\phi^{n-\frac12}.
\]
Without losing generality, we consider
periodic boundary conditions with a proper initial data $\phi^0(\mathbf{x})$.
Our analysis can be extended in a straightforward way to the
fully discrete numerical schemes with some appropriate spatial discretization
preserving the discrete integration-by-parts formulas,
such as the Fourier pseudo-spectral method
\cite{ChengWangWise:2016WeaklyEnergy,ChengWangWise:2019SPFC}
used in our experiments.

A major hindrance in establishing the discrete energy stability
for the time-fractional phase field models is
the positive definiteness of the following real quadratic form
with respect to the convolution kernels $a_{n-j}^{(n)}$ (in a general sense)
arising from variable-step time approximations
\begin{align}\label{equality:nonuniform L1 real quadratic}
\lan{2}\sum_{k=1}^nw_k\sum_{j=1}^ka_{k-j}^{(k)}w_j\quad
\text{for any nonzero sequence $\{w_1,w_2,\cdots,w_n\}$}.
\end{align}
On the uniform time mesh with $a_{j}^{(n)}=a_{j}$,
L\'opez-Marcos \cite[Proposition 5.2]{Lopez:1990ADifferenceScheme}
gave some sufficient but algebraic conditions to check the desired property,
\begin{align}\label{condition:unifrom real quadratic}
a_j\ge 0,\quad a_{j-1}\ge a_j\quad\text{and}\quad a_{j-1}-a_j\ge a_j-a_{j+1}.
\end{align}
The positive, decreasing and convex criteria have been widely used
to establish the stability and convergence results
for integro-differential and time-fractional differential problem,
such as \cite[Lemma 2.6]{JiLiaoZhang:2020Simple} for a discrete (global) energy law
for the time-fractional Allen-Cahn model.
\lan{Very recently, Karaa \cite{Karaa:2021UniformPositivityL1Formula} presented some criteria ensuring the positivity of the real quadratic form  \eqref{equality:nonuniform L1 real quadratic} for
some commonly used numerical methods, including the convolution quadrature method and L1 formula.
}

\lan{It is worthwhile noting that the criterion \eqref{condition:unifrom real quadratic} 
may fail to verify the desired positive definiteness of
the real quadratic form  \eqref{equality:nonuniform L1 real quadratic} for the
variable-step time-stepping methods,
such as the nonuniform L1 method \cite[Remark 1]{JiLiaoZhang:2020Simple}.
Also, the  technique of completely monotone sequence in \cite{Karaa:2021UniformPositivityL1Formula} may not be applied to the nonuniform case directly.
} 
Recently, Liao et. al. \cite{LiaoTangZhou:2020Positive}
proposed another class of sufficient but easy-to-check  conditions
(for general discrete kernels)
\begin{align}\label{condition:nonunifrom real quadratic}
a_{j-1}^{(n)}\ge a_{j}^{(n)}>0,\quad
a_{j-1}^{(n-1)}>a_{j}^{(n)}\quad\text{and}\quad
a_{j-1}^{(n-1)}a_{j+1}^{(n)}\ge a_{j}^{(n-1)}a_{j}^{(n)}.
\end{align}
The main theorem of \cite[Theorem 1.1]{LiaoTangZhou:2020Positive} ensures that
the positive, decreasing and convex criteria \eqref{condition:nonunifrom real quadratic}
are sufficient for the positive definiteness
of associated quadratic form resulting from a general class of discrete convolution approximations.
By a careful verification of the sufficient condition
\eqref{condition:nonunifrom real quadratic},
the positive definiteness of the nonuniform L1 kernels
$a_{n-k}^{(n)}$ was verified in
\cite[Proposition 4.1]{LiaoTangZhou:2020Positive}.
As a direct application, the stabilized semi-implicit
scheme was shown to preserve the global energy stability
on arbitrary meshes, see \cite[Proposition 4.2]{LiaoTangZhou:2020Positive}.

\begin{figure}[htb!]
\centering
\includegraphics[width=4.4in]{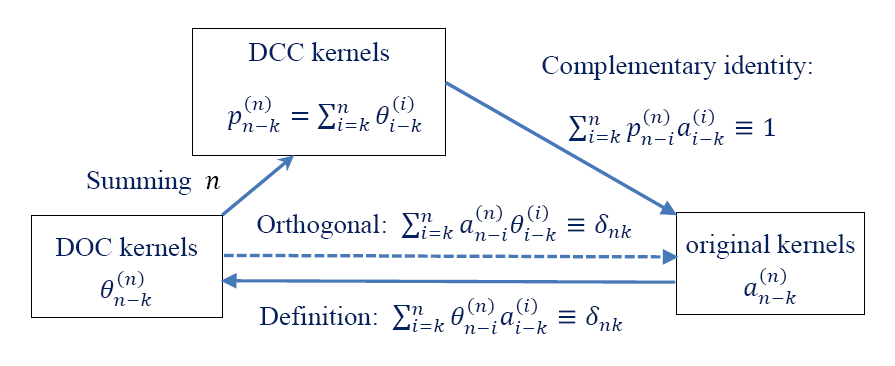}
\caption{The relationship diagram of DOC, DCC and original kernels.}
\label{fig: DOC and DCC relation}
\end{figure}

However, the positive
definiteness property in \cite[Proposition 4.1]{LiaoTangZhou:2020Positive}
\lan{may not be sufficient} to build up a discrete energy
stability for some fully implicit methods,
such as the implicit L1 scheme (2.11) in \cite{JiLiaoZhang:2020Simple}.
By means of some new discrete tools,
the recent work \cite{LiaoZhuWang:2020TFAC}
filled this gap for the time-fractional Allen-Cahn model
by demonstrating that the implicit L1 scheme
preserving the discrete variational energy law
\eqref{cont:TFCH variational energy}
on arbitrary time meshes.
One of the key tools is the discrete orthogonal convolution
(DOC) kernels defined by the following recursive procedure
\begin{align}\label{def:L1 DOC kernels}
\theta_{0}^{(n)}:=\frac{1}{a_{0}^{(n)}}
\quad \mathrm{and} \quad
\theta_{n-k}^{(n)}:=-\frac{1}{a_{0}^{(k)}}
\sum_{j=k+1}^n {\theta_{n-j}^{(n)}}a_{j-k}^{(j)}
\quad \text{for $1\le k\le n-1$}.
\end{align}
It is easy to check that the DOC kernels $\theta_{n-k}^{(n)}$ satisfy
the discrete orthogonal identity
\begin{align}\label{equality:orthogonal identity}
\sum_{j=k}^n\theta_{n-j}^{(n)}a_{j-k}^{(j)}\equiv \delta_{nk}
\quad\text{for $1\le k\le n$,}
\end{align}
where $\delta_{nk}$ is the Kronecker delta symbol.
Furthermore, another useful discrete analysis tool is the so-called
discrete complementary convolution (DCC)
kernels introduced by means of the DOC kernels $\theta_{n-k}^{(n)}$,
see \cite[Subsection 2.2]{LiaoTangZhou:2020Positive},
\begin{align}\label{def:L1 DCC kernels}
p_{n-k}^{(n)}:=\sum_{j=k}^n\theta_{j-k}^{(j)}\quad
\text{such that}\quad \sum_{j=k}^np_{n-j}^{(n)}a_{j-k}^{(j)}\equiv 1
\quad \text{for $1\le k\le n$}.
\end{align}
The interplay relationship of the mentioned DOC and DCC kernels together with the original kernels
is summarized in Figure \ref{fig: DOC and DCC relation}, also see \cite[Figure 1]{LiaoTangZhou:2020Positive}.

In this paper, we investigate the positive definiteness
of the L1-type kernels from the L1-type approximations
\eqref{scheme:L1 formula}, \eqref{scheme:L1h formula} and \eqref{scheme:L1a formula}
of the Caputo derivative, and explore the energy stability of the associated numerical methods
\eqref{scheme:fully implicit L1 TFCH}-\eqref{scheme:fully implicit L1a TFCH}.
Some novel discrete convolution inequalities with respect to
the L1 and L1$_{\mathrm{h}}$ kernels are established in Theorems \ref{thm:L1 kernel positive definite}
and \ref{thm:L1h kernel positive definite}, respectively.
\lan{Let $\lambda_{\min}$ be the minimum eigenvalue of the real quadratic form (matrix) involving the L1 kernels ($\lambda_{\min}^{(\mathrm{h})}$ and $\lambda_{\min}^{(\mathrm{a})}$ of L1$_{\mathrm{h}}$ kernels and L1$_{\mathrm{a}}$ kernels are defined similarly).}
Subsection 2.1 obtains the certain lower bound of the minimum eigenvalue.
To the best of our knowledge, such estimate of the minimum eigenvalue on arbitrary time meshes
is considered at the first time.
The positive definiteness of the L1$_{\mathrm{h}}$ kernels \eqref{def:L1h formula kernel}
is also verified in subsection 2.2 although \lan{the first two kernels lose their monotonicity},
see Table \ref{table:comparison L1-type formula} which collects some related properties of
the underlying discrete kernels. However, the positive definite property
of the averaged L1$_{\mathrm{a}}$ kernels \eqref{def:L1a formula kernel}
is still undetermined, see more details in subsection 2.3.

\begin{table}[htb!]
\begin{center}
\caption{Kernels properties of L1-type formulas on nonuniform time meshes.}
\label{table:comparison L1-type formula} \vspace*{0.3pt}
\def\temptablewidth{0.9\textwidth}
{\rule{\temptablewidth}{1pt}}
\begin{spacing}{1.35}
\begin{tabular*}{\temptablewidth}{@{\extracolsep{\fill}}l|ccc}
  Theoretical properties &L1 kernels \eqref{def:L1 formula kernel}   &L1$_\mathrm{h}$ kernels
  \eqref{def:L1h formula kernel}   &L1$_\mathrm{a}$ kernels \eqref{def:L1a formula kernel}\\
  \midrule
  Positivity &$a_j^{(n)}>0$  &$\Lh{j}{n}>0$ &$\La{j}{n}>0$\\[4pt]
  Monotonicity &strict decreasing  &$\Lh{0}{n}\ngeq\Lh{1}{n}$
  &$\La{0}{n}\ngeq\La{1}{n}$ \\
  \tabincell{l}{Positive definiteness \\(with eigenvalues $\lambda$)}
  &$\lambda_{\min}\ge \min\limits_{1\le k\le n}a_{0}^{(k)}$  &$\lambda_{\min}^{(\mathrm{h})}>0$&\textsf{undetermined}\\
\end{tabular*}
\end{spacing}
{\rule{\temptablewidth}{0.5pt}}
\end{center}
\end{table}	

Then we show in Theorems \ref{thm:L1h energy law TFCH}
and \ref{thm:L1h energy law TFCH} that
the proposed numerical schemes
\eqref{scheme:fully implicit L1 TFCH}-\eqref{scheme:fully implicit L1h TFCH}
preserve the variational energy dissipation law \eqref{cont:TFCH variational energy}
on arbitrary time meshes.
These discrete energy laws
are shown to be asymptotically compatible with the classical energy dissipation laws
as the fractional order $\alpha\rightarrow1$, see Remarks
\ref{remark:L1 energy compatibility}-\ref{remark:L1h energy compatibility}.


The rest of this paper is organized as follows.
In section 2, the positive definiteness of
the suggested L1-type formulas
is investigated.
Section 3 establishes the discrete energy dissipation laws of the  L1-type time-stepping schemes.
Numerical examples are presented in section 4 to confirm our theoretical findings.

\section{Positive definiteness of L1-type formulas}
\setcounter{equation}{0}

\subsection{Positive definiteness of L1 kernels}
It follows from \cite[Proposition 4.1]{LiaoTangZhou:2020Positive}
that the discrete L1 kernels $a_{n-k}^{(n)}$
satisfy the criteria \eqref{condition:nonunifrom real quadratic}
so that the following result holds according to
\cite[Theorem 1.1]{LiaoTangZhou:2020Positive}.

\begin{lemma}\label{lem:old L1 kernel positive definiteness}
The L1 kernels $a_{j}^{(n)}$ in \eqref{def:L1 formula kernel}  are positive definite in the sense that
\begin{align}\label{equality:L1 kernel real quadratic}
\lan{2}\sum_{k=1}^nw_k\sum_{j=1}^ka_{k-j}^{(k)}w_j>0\quad
\text{for any nonzero sequence $\{w_1,w_2,\cdots,w_n\}$}.
\end{align}
\end{lemma}

In what follows, we improve Lemma \ref{lem:old L1 kernel positive definiteness}
by presenting a lower bound $\sigma_{\text{L1}}>0$ for
the minimum eigenvalue of the associated quadratic form in the sense that
\begin{align}\label{inequality:lower bound L1 kernel}
\lan{2}\sum_{k=1}^nw_k\sum_{j=1}^ka_{k-j}^{(k)}w_j
\ge \sigma_{\text{L1}}\sum_{k=1}^nw_k^2.
\end{align}
To this end, we first list some properties of
the DOC kernels $\theta_{n-k}^{(n)}$ defined in
\eqref{def:L1 DOC kernels},
and the DCC kernels $p_{n-k}^{(n)}$ defined in \eqref{def:L1 DCC kernels},
see \cite[Lemmas 2.3 and 2.5]{LiaoTangZhou:2020Positive}.
\begin{lemma} \label{lem:L1 DOC DCC kernel property}
For any $n\ge 2$, the DOC kernels $\theta_{n-k}^{(n)}$ in \eqref{def:L1 DOC kernels} satisfy
\[
\theta_{0}^{(n)}>0\quad \text{and}\quad \theta_{n-k}^{(n)}<0
\quad \text{for $1\le k\le n-1$};\quad
\text{but} \quad\sum_{k=1}^n\theta_{n-k}^{(n)}>0;
\]
and the DCC kernels $p_{n-k}^{(n)}$ in \eqref{def:L1 DCC kernels} satisfy
\[p_{n-k}^{(n)}=\sum_{j=k}^n\theta_{j-k}^{(j)}\ge 0
\quad \text{for $1\le k\le n$.}
\]
\end{lemma}

According to the definition \eqref{def:L1 DCC kernels},
one can find that the DOC kernels $\theta_{n-k}^{(n)}$ and the DCC
kernels $p_{n-k}^{(n)}$ satisfy the following relationship
\begin{align}\label{def:L1 DCC-DOC relationship}
\theta_{0}^{(n)}=p_{0}^{(n)}\quad \text{and}\quad
\theta_{n-k}^{(n)}=p_{n-k}^{(n)}-p_{n-1-k}^{(n-1)}
\quad \text{for $1\le k\le n-1$}.
\end{align}
Then we have the following  positive definiteness result
for the DOC kernels $\theta_{n-k}^{(n)}$.
\begin{lemma}\label{lem:L1 DOC kernel positive definite}
For any real vector sequence $\{w_k\}_{k=1}^n$, it holds that
\begin{align*}
2w_k\sum_{j=1}^k\theta_{k-j}^{(k)}w_j
\ge \sum_{j=1}^kp_{k-j}^{(k)}w_j^2-\sum_{j=1}^{k-1}p_{k-1-j}^{(k-1)}w_j^2
+\frac{1}{\theta_0^{(k)}}\braB{\sum_{j=1}^k\theta_{k-j}^{(k)}w_j}^2
\quad \text{for $ k\ge 1$,}
\end{align*}
so that the DOC kernels $\theta_{n-k}^{(n)}$ are positive definite in the sense that
\begin{align*}
2\sum_{k=1}^nw_k\sum_{j=1}^k\theta_{k-j}^{(k)}w_j
\ge \sum_{k=1}^np_{n-k}^{(n)}w_k^2
+\sum_{k=1}^n\frac{1}{\theta_0^{(k)}}\braB{\sum_{j=1}^k\theta_{k-j}^{(k)}w_j}^2
>0\quad \text{if $w_k\not\equiv0$.}
\end{align*}
\end{lemma}
\begin{proof}
\lan{The first inequality can be verified by the proof of
\cite[Lemma 2.4]{LiaoZhuWang:2020TFAC}.
Summing up this inequality from $k=1$ to $n$,
one gets the claimed second inequality and the proof is completed.}
\end{proof}

\begin{theorem}\label{thm:L1 kernel positive definite}
For any real vector sequence $\{w_k\}_{k=1}^n$, it holds that
\begin{align*}
2w_k\sum_{j=1}^ka_{k-j}^{(k)}w_j
\ge a_0^{(k)}w_k^2
+\sum_{j=1}^kp_{k-j}^{(k)}\braB{\sum_{\ell=1}^ja_{j-\ell}^{(j)}w_\ell}^2
-\sum_{j=1}^{k-1}p_{k-1-j}^{(k-1)}\braB{\sum_{\ell=1}^ja_{j-\ell}^{(j)}w_\ell}^2
\end{align*}
for $ k\ge 1$, so that the L1 kernels $a_{n-k}^{(n)}$ in \eqref{def:L1 formula kernel} are
positive definite in the sense that
\begin{align*}
2\sum_{k=1}^nw_k\sum_{j=1}^ka_{k-j}^{(k)}w_j
\ge \sum_{k=1}^na_{0}^{(k)}w_k^2
+\sum_{k=1}^np_{n-k}^{(n)}\braB{\sum_{j=1}^ka_{k-j}^{(k)}w_j}^2
>0\quad \text{if $w_k\not\equiv0$.}
\end{align*}
\end{theorem}
\begin{proof}
For any fixed index $n\ge1$ and any real vector sequence $\{w_k\}_{k=1}^n$,
we introduce an auxiliary sequence $\{v_k\}_{k=1}^n$ by means of the
original kernels $a_{j-k}^{(j)}$ as follows
\begin{align}\label{lemProof:L1 kernel auxiliary sequence}
v_j:=\sum_{k=1}^ja_{j-k}^{(j)}w_k
\quad\text{for $1\le j\le n$.}
\end{align}
Multiplying both sides of the equality \eqref{lemProof:L1 kernel auxiliary sequence}
by $\theta_{j-k}^{(j)}$ and summing up from $k=1$ to $j$ give
\begin{align}\label{lemProof:relation wn-vn}
\sum_{k=1}^j\theta_{j-k}^{(j)}v_k
=\sum_{k=1}^j\theta_{j-k}^{(j)}\sum_{\ell=1}^ka_{k-\ell}^{(k)}w_\ell
=\sum_{\ell=1}^jw_\ell\sum_{k=\ell}^j\theta_{j-k}^{(j)}a_{k-\ell}^{(k)}
=w_j\quad\text{for $1\le j\le n$,}
\end{align}
where the discrete orthogonal identity \eqref{equality:orthogonal identity}
has been used in the last step. According to the first inequality in
Lemma \ref{lem:L1 DOC kernel positive definite} and
the fact $a_0^{(k)}=1/\theta_0^{(k)}$ from \eqref{def:L1 DOC kernels}, one has
\begin{align*}
2v_k\sum_{j=1}^k\theta_{k-j}^{(k)}v_j
\ge \sum_{j=1}^kp_{k-j}^{(k)}v_j^2-\sum_{j=1}^{k-1}p_{k-1-j}^{(k-1)}v_j^2
+a_0^{(k)}\braB{\sum_{j=1}^k\theta_{k-j}^{(k)}v_j}^2
\quad \text{for $ k\ge 1$,}
\end{align*}
which directly leads to the claimed first inequality  by the above relationships
\eqref{lemProof:L1 kernel auxiliary sequence}-\eqref{lemProof:relation wn-vn}.
\end{proof}

\lan{Theorem \ref{thm:L1 kernel positive definite} expresses the discrete convolution structure of the nonuniform L1 formula \eqref{scheme:L1 formula} with the original convolution kernels 
	$a_{n-k}^{(n)}$ rather than the corresponding DOC kernels $\theta_{n-k}^{(n)}$. This form will be heuristic in treating other numerical Caputo derivatives, especially when the associated discrete kernels lose the monotonicity, see the L1$_\mathrm{h}$ kernels in next subsection. As a byproduct, this form updates Lemma \ref{lem:old L1 kernel positive definiteness} by presenting a lower bound
$$\sigma_{\text{L1}}:=\min_{1\le k\le n} a_{0}^{(k)}$$ for
the minimum eigenvalue $\lambda_{\min}$ of the associated quadratic form.
Table \ref{table:min eigenvalue L1 kernel}
tabulates the low bound $\sigma_{\text{L1}}$ and the minimum eigenvalue $\lambda_{\min}$
on the random time meshes with $T=1$
for three different fractional orders $\alpha$. 
As observed, $\sigma_{\text{L1}}$ is a delicate estimate of $\lambda_{\min}$,
especially when the fractional order $\alpha$ is small. Table \ref{table:min eigenvalue L1 kernel on garded meshes} also lists the comparsions on the graded mesh $t_k=T(k/N)^\gamma$, which was applied frequently in resolving the initial singularity  \cite{Stynes:2017ErrorAnalysis,Kopteva:2019Error,LiaoLiZhang:2018SharpError}.}

\begin{table}[htb!]
\begin{center}
\caption{Comparisons of the bound
 $\sigma_{\text{L1}}$ and $\lambda_{\min}$ on random time meshes.}
\label{table:min eigenvalue L1 kernel} \vspace*{0.3pt}
\def\temptablewidth{0.6\textwidth}
{\rule{\temptablewidth}{0.5pt}}
\begin{spacing}{1.25}
\begin{tabular*}{\temptablewidth}{@{\extracolsep{\fill}}c|cc|cc|cc}
  \multirow{2}{*}{$N$} &\multicolumn{2}{c|}{$\alpha=0.1$}
  &\multicolumn{2}{c|}{$\alpha=0.5$} &\multicolumn{2}{c}{$\alpha=0.9$}\\
  \cline{2-3} \cline{4-5} \cline{6-7}
  &$\sigma_{\text{L1}}$ &$\lambda_{\min}$
  &$\sigma_{\text{L1}}$ &$\lambda_{\min}$
  &$\sigma_{\text{L1}}$ &$\lambda_{\min}$\\
  \midrule
  100	 &1.53 &1.72	   &7.73  &12.37    &33.57  &65.85   \\
  200	 &1.65 &1.84	   &11.19 &17.91    &65.35  &128.74  \\
  400	 &1.76 &1.97	   &15.83 &24.99    &121.98 &239.22  \\
\end{tabular*}
\end{spacing}
{\rule{\temptablewidth}{0.5pt}}
\end{center}
\end{table}

\begin{table}[htb!]
	\begin{center}
		\caption{\lan{Comparisons  of the bound
				$\sigma_{\text{L1}}$ and $\lambda_{\min}$ for $\alpha=0.5$ on graded time meshes.}}
		\label{table:min eigenvalue L1 kernel on garded meshes} \vspace*{0.3pt}
		\def\temptablewidth{0.6\textwidth}
		{\rule{\temptablewidth}{0.5pt}}
		\begin{spacing}{1.25}
			\begin{tabular*}{\temptablewidth}{@{\extracolsep{\fill}}c|cc|cc|cc}
				\multirow{2}{*}{$N$} &\multicolumn{2}{c|}{$\gamma=1$}
				&\multicolumn{2}{c|}{$\gamma=2$} &\multicolumn{2}{c}{$\gamma=4$}\\
				\cline{2-3} \cline{4-5} \cline{6-7}
				&$\sigma_{\text{L1}}$ &$\lambda_{\min}$
				&$\sigma_{\text{L1}}$ &$\lambda_{\min}$
				&$\sigma_{\text{L1}}$ &$\lambda_{\min}$\\
				\midrule
				100	 &11.28	 &17.16   &8.00	  &12.36  &5.68	  &8.92  \\
				200	 &15.96  &24.26	  &11.30  &17.36  &8.01	  &12.43  \\
				400	 &22.57	 &34.31   &15.97  &24.44  &11.30  &17.42  \\
			\end{tabular*}
		\end{spacing}
		{\rule{\temptablewidth}{0.5pt}}
	\end{center}
\end{table}

The sharpness of the bound $\sigma_{\text{L1}}$
can be also seen by comparing with a previous result in \cite[Lemma 3.1]{TangYuZhou:2018TimeFractionalPhase} on the uniform grid. In this simple case, Theorem \ref{thm:L1 kernel positive definite} shows $\sigma_{\text{L1}}
=\frac{1}{\Gamma(2-\alpha){\tau}^{\alpha}}$, while 
Tang et al \cite[Lemma 3.1]{TangYuZhou:2018TimeFractionalPhase} gave the following bound
$$\sigma_{*}:=\frac{1}{\tau^{\alpha}\Gamma(1-\alpha)}\brab{\frac{2}{n+1}}^\alpha.$$
The current estimate $\sigma_{\text{L1}}$
is sharper than $\sigma_{*}$, that is, 
$$\frac{\sigma_{\text{L1}}}{\sigma_{*}}
=\frac{2^{-\alpha}}{1-\alpha}\brat{n+1}^\alpha>\brat{n+1}^\alpha,$$
due to the fact
$2^{-\alpha}>1-\alpha$ for $\alpha\in(0,1)$.
\lan{Very recently, Karaa\cite[Lemma 3.4]{Karaa:2021UniformPositivityL1Formula}  gave a new bound $\sigma_{\star}$ on the uniform mesh, that is, 
	\begin{align*}
	\sigma_{\star}:=\frac{2 \mathrm{Li}_{\alpha-1}(-1)}{(-1) \Gamma(2-\alpha)\tau^{\alpha}},
	\end{align*}
where the polylogarithm function $\mathrm{Li}_{\beta}(z):=\sum_{j=1}^{\infty}\frac{z^j}{j^{\beta}}$, which is well defined for $|z|<1$ and can be analytically
extended to the split domain $\mathbb{C}\backslash[1,\infty)$.
Actually, $\sigma_{\star}\ge \sigma_{\text{L1}}$ on the uniform mesh due to the fact 
$(-1)\mathrm{Li}_{\alpha-1}(-1)\ge 1/4$.
It seems that the result of Theorem \ref{thm:L1 kernel positive definite} still has a lot of room for improvement.}

\subsection{Positive definiteness of L1$_\mathrm{h}$ kernels}

The L1$_\mathrm{h}$ kernels in \eqref{def:L1h formula kernel} are different from the L1 kernels
due to the lack of the monotonicity.
Simple calculations show that the first
two discrete L1$_\mathrm{h}$ kernels defined in \eqref{def:L1h formula kernel} satisfy
\begin{align*}
\Lh{0}{n}-\Lh{1}{n}
=\Lh{0}{n}\kbraB{1+r_n-r_n^\alpha\brab{r_n+2}^{1-\alpha}}.
\end{align*}
It is evident that  $\Lh{0}{n}> \Lh{1}{n}$ as $\alpha\rightarrow 1$
and $\Lh{0}{n}< \Lh{1}{n}$ as $\alpha\rightarrow 0,$ i.e.,
the monotonously decreasing
of the L1$_{\mathrm{h}}$ kernels $\Lh{j}{n}$
loses for some fractional orders $\alpha\in(0,1)$.
So  the sufficient criterion \eqref{condition:nonunifrom real quadratic}
or \cite[Theorem 1.1]{LiaoTangZhou:2020Positive}
can not be directly applied to confirm the positive definiteness
of the discrete L1$_\mathrm{h}$ kernels $a_{j}^{(\mathrm{h},n)}$.

Further observations suggest that the desired  monotonicity property
can be attained by doubling the first kernel $\Lh{0}{n}$. Actually,
the inequality $1+\beta z>(1+z)^\beta$ holds for any $\beta\in(0,1)$
with respect to $z>0$ and then
\begin{align}\label{inequality:first two L1h kernel}
2\Lh{0}{n}-\Lh{1}{n}
&=\Lh{0}{n}\kbrab{2+r_n-r_n^\alpha\brab{r_n+2}^{1-\alpha}}\nonumber\\
&>\Lh{0}{n}r_n\kbra{1+2(1-\alpha)/r_n-(1+2/r_n)^{1-\alpha}}>0.
\end{align}
Accordingly, we introduce the following auxiliary kernels
via the original L1$_\mathrm{h}$ kernels $\Lh{j}{n}$,
\begin{align}\label{def:auxiliary L1h kernel}
\aLh{a}{0}{n}:=2\Lh{0}{n}\quad
\text{and}\quad \aLh{a}{j}{n}:=\Lh{j}{n}
 \quad\text{for $1\le j\le n-1$.}
\end{align}

\begin{lemma}\label{lem:L1h auxiliary  kernel positive definiteness}
For any fixed $n\ge 2$, the auxiliary L1$_{\mathrm{h}}$ kernels $\aLh{a}{j}{n}$
in \eqref{def:auxiliary L1h kernel} satisfy
\begin{align*}
\aLh{a}{j-1}{n}\ge\aLh{a}{j}{n}>0,\quad
\aLh{a}{j-1}{n-1}>\aLh{a}{j}{n},\quad
\aLh{a}{j-1}{n-1}\aLh{a}{j+1}{n}
>\aLh{a}{j}{n-1}\aLh{a}{j}{n}.
\end{align*}
Then the auxiliary L1$_{\mathrm{h}}$ kernels $\aLh{a}{j}{n}$
are positive definite according to \cite[Theorem 1.1]{LiaoTangZhou:2020Positive}.
\end{lemma}
\begin{proof}
The positivity and the monotonous  decreasing of
$\aLh{a}{j}{n}$
follow from  the definition \eqref{def:auxiliary L1h kernel}
and the integral mean value theorem immediately.
Then, we consider the following sequence
\[
\psi_{n-k}^{(\mathrm{h},n)}:=\frac{\aLh{a}{n-k}{n}}{\aLh{a}{n-1-k}{n-1}}
=\frac{b_{n,k}(1)}{b_{n-1,k}(1)}
\quad \text{for $1\le k\le n-1$,}
\]
where the auxiliary functions $b_{n,k}(\xi)$ are defined by
\[
b_{n,k}(\xi):=\frac{1}{\tau_k}\int_{t_{k-1}}^{t_{k-1}+\tau_k\xi}
\omega_{1-\alpha}(t_{n-\frac12}-s)\zd{s}\quad
\text{for $1\le k\le n-1$},
\]
and
\[
b_{n,n}(\xi)
:=\frac{2}{\tau_n}\int_{t_{n-1}}^{t_{n-1}+\frac{\tau_n}{2}\xi}
\omega_{1-\alpha}(t_{n-\frac12}-s)\zd{s}\quad \text{for $k=n$}.
\]
Differentiating the functions $b_{n,k}(\xi)$ yields
$b^{\prime}_{n,n}(\xi)=2^\alpha\omega_{1-\alpha}(\tau_n-\tau_n\xi)$
and
\[
b^{\prime}_{n,k}(\xi)=\omega_{1-\alpha}(t_{n-\frac12}-t_{k-1}-\tau_k\xi)
\quad\text{for $1\le k\le n-1$}.
\]
Then the  Cauchy mean value theorem shows  that
there exists some $\xi_{n-1}\in(0,1)$ such that
\begin{align*}
\psi_{1}^{(\mathrm{h},n)}
&=\frac{b_{n,n-1}(1)-b_{n,n-1}(0)}{b_{n-1,n-1}(1)-b_{n-1,n-1}(0)}
=\frac{b_{n,n-1}^\prime(\xi_{n-1})}
{b_{n-1,n-1}^\prime(\xi_{n-1})}\\
&=\frac{1}{2^\alpha}\braB{\frac{\tau_{n-1}-\tau_{n-1}\xi_{n-1}}
{\tau_n/2+\tau_{n-1}-\tau_{n-1}\xi_{n-1}}}^\alpha
<\braB{\frac{\tau_{n-1}}{\tau_n+2\tau_{n-1}}}^\alpha.
\end{align*}
Here, we use the fact that the function
$y=(A-z)/(B-z)$ is monotonically decreasing with respect
to the variable $z$ if the two parameters $A<B$.
Analogously, one can follow the above proof or
\cite[Proposition 4.1]{LiaoTangZhou:2020Positive}
to derive
\[
\braB{\frac{t_{n-\frac32}-t_k}{t_{n-\frac12}-t_k}}^\alpha
<\psi_{n-k}^{(\mathrm{h},n)}
<\braB{\frac{t_{n-\frac32}-t_{k-1}}{t_{n-\frac12}-t_{k-1}}}^\alpha
\quad \text{for $1\le k\le n-2$.}
\]
Thus we have the following inequalities
\[
\psi_{1}^{(\mathrm{h},n)}<\braB{\frac{\tau_{n-1}}{\tau_n+2\tau_{n-1}}}^\alpha
<\psi_{2}^{(\mathrm{h},n)}<\psi_{3}^{(\mathrm{h},n)}
<\cdots<\psi_{n-1}^{(\mathrm{h},n)}
<1\quad \text{for $n\ge 2$.}
\]
They lead to the remainder properties (the last two classes of inequalities) for the auxiliary L1$_{\mathrm{h}}$ kernels $\aLh{a}{j}{n}$. Then
\cite[Theorem 1.1]{LiaoTangZhou:2020Positive} completes the proof.
\end{proof}

The above results show that the auxiliary L1$_{\mathrm{h}}$ kernels $\aLh{a}{n-k}{n}$
satisfy the criterion \eqref{condition:nonunifrom real quadratic}.
It is reasonable to define the associated DOC kernels $\aLh{\uptheta}{n-k}{n}$ as follows,
\begin{align}\label{def:auxiliary L1h DOC kernels}
\aLh{\uptheta}{0}{n}:=\frac{1}{\aLh{a}{0}{n}}
\quad \mathrm{and} \quad
\aLh{\uptheta}{n-k}{n}:=-\frac{1}{\aLh{a}{0}{k}}
\sum_{j=k+1}^n \aLh{\uptheta}{n-j}{n}\aLh{a}{j-k}{j}
\quad \text{for $1\le k\le n-1$}.
\end{align}
Also, we can define the corresponding DCC kernels $\aLh{p}{n-k}{n}$ by
\begin{align}\label{def:auxiliary L1h DCC kernels}
\aLh{p}{n-k}{n}=\sum_{j=k}^n\aLh{\uptheta}{j-k}{j}
\quad \text{for $1\le k\le n$}.
\end{align}
They are well-defined and satisfy the following results according to
\cite[Lemmas 2.3 and 2.5]{LiaoTangZhou:2020Positive}.

\begin{lemma} \label{lem:L1h DOC DCC kernel property}
For any $n\ge 2$, the DOC kernels $\aLh{\uptheta}{n-k}{n}$ in
\eqref{def:auxiliary L1h DOC kernels} satisfy
\[
\aLh{\uptheta}{0}{n}>0\quad
\text{and}\quad \aLh{\uptheta}{n-k}{n}<0
\quad \text{for $1\le k\le n-1$}\quad
\text{but} \quad\sum_{k=1}^n\aLh{\uptheta}{n-k}{n}>0;
\]
and the DCC kernels $\aLh{p}{n-k}{n}$ in
\eqref{def:auxiliary L1h DCC kernels} satisfy
\[\aLh{p}{n-k}{n}=\sum_{j=k}^n\aLh{\uptheta}{j-k}{j}\ge 0
\quad \text{for $1\le k\le n$}.
\]
\end{lemma}

By using Lemma \ref{lem:L1h DOC DCC kernel property},
one can follow the proof of Lemma \ref{lem:L1 DOC kernel positive definite}
to prove the following result.
\lan{
\begin{lemma}\label{lem:L1h auxiliary DOC kernel positive definite}
For any real vector sequence $\{w_k\}_{k=1}^n$, it holds that
\begin{align*}
2w_k\sum_{j=1}^k\aLh{\uptheta}{k-j}{k}w_j
\ge \sum_{j=1}^k\aLh{p}{k-j}{k}w_j^2
-\sum_{j=1}^{k-1}\aLh{p}{k-1-j}{k-1}w_j^2
+\frac{1}{\aLh{\uptheta}{0}{k}}
\braB{\sum_{j=1}^k\aLh{\uptheta}{k-j}{k}w_j}^2
\quad \text{for $ k\ge 1$.}
\end{align*}
\end{lemma}
}

We are in a position to verify the positive definiteness of the
L1$_\mathrm{h}$ kernels $\Lh{n-k}{n}$ in \eqref{def:L1h formula kernel}.
\begin{theorem}\label{thm:L1h kernel positive definite}
For any real vector sequence $\{w_k\}_{k=1}^n$, it holds that
\begin{align*}
2w_k\sum_{j=1}^k\Lh{k-j}{k}w_j
\ge \sum_{j=1}^k\aLh{p}{k-j}{k}\braB{\sum_{\ell=1}^j\aLh{a}{j-\ell}{j}w_\ell}^2
-\sum_{j=1}^{k-1}\aLh{p}{k-1-j}{k-1}\braB{\sum_{\ell=1}^j\aLh{a}{j-\ell}{j}w_\ell}^2
\quad \text{for $ k\ge 1$,}
\end{align*}
so that the L1$_\mathrm{h}$ kernels $\Lh{n-k}{n}$ in \eqref{def:L1h formula kernel}
are positive definite in the sense that
\begin{align*}
2\sum_{k=1}^nw_k\sum_{j=1}^k\Lh{k-j}{k}w_j
\ge\sum_{k=1}^n\aLh{p}{n-k}{n}\braB{\sum_{j=1}^k\aLh{a}{k-j}{k}w_j}^2
>0\quad \text{for $ n\ge 1$ if $w_k\not\equiv0$.}
\end{align*}
\end{theorem}
\begin{proof}
By using Lemma \ref{lem:L1h auxiliary DOC kernel positive definite},
one can follow the proof of Theorem \ref{thm:L1 kernel positive definite} to obtain
\begin{align*}
2w_k\sum_{j=1}^k\aLh{a}{k-j}{k}w_j
\ge \aLh{a}{0}{k}w_k^2
+\sum_{j=1}^k\aLh{p}{k-j}{k}
\braB{\sum_{\ell=1}^j\aLh{a}{j-\ell}{j}w_\ell}^2
-\sum_{j=1}^{k-1}\aLh{p}{k-1-j}{k-1}
\braB{\sum_{\ell=1}^j\aLh{a}{j-\ell}{j}w_\ell}^2.
\end{align*}
Then the definition \eqref{def:auxiliary L1h kernel}
of $\aLh{a}{k-j}{k}$ implies
\begin{align*}
2w_k\sum_{j=1}^k\Lh{k-j}{k}w_j
\ge\sum_{j=1}^k\aLh{p}{k-j}{k}
\braB{\sum_{\ell=1}^j\aLh{a}{j-\ell}{j}w_\ell}^2
-\sum_{j=1}^{k-1}\aLh{p}{k-1-j}{k-1}
\braB{\sum_{\ell=1}^j\aLh{a}{j-\ell}{j}w_\ell}^2.
\end{align*}
Summing up this inequality from $k=1$ to $n$ yields the claimed result.
\end{proof}
%

\begin{table}[htb!]
\begin{center}
\caption{The minimum eigenvalue $\lambda_{\min}^{(\mathrm{h})}$  on random time meshes.}
\label{table:L1h kernel verification} \vspace*{0.3pt}
\def\temptablewidth{0.6\textwidth}
{\rule{\temptablewidth}{0.5pt}}
\begin{spacing}{1.25}
\begin{tabular*}{\temptablewidth}{@{\extracolsep{\fill}}c|ccc}
$N$    &$\alpha=0.1$  &$\alpha=0.5$ &$\alpha=0.9$\\
  \midrule
  100    &0.20	        &6.64         &60.69 \\	 	
  200    &0.21	        &9.48         &119.07 \\	 	
  400    &0.23	        &13.09        &219.60\\
\end{tabular*}
\end{spacing}
{\rule{\temptablewidth}{0.5pt}}
\end{center}
\end{table}

\begin{table}[htb!]
	\begin{center}
		\caption{\lan{The minimum eigenvalue $\lambda_{\min}^{(\mathrm{h})}$  on graded time meshes.}}
		\label{table:L1h kernel verification on graded meshes} \vspace*{0.3pt}
		\def\temptablewidth{0.65\textwidth}
		{\rule{\temptablewidth}{0.5pt}}
		\begin{spacing}{1.25}
			\begin{tabular*}{\temptablewidth}{@{\extracolsep{\fill}}c|cc|cc|cc}
				\multirow{2}{*}{$N$} &\multicolumn{2}{c|}{$\alpha=0.1$}
				&\multicolumn{2}{c|}{$\alpha=0.5$} &\multicolumn{2}{c}{$\alpha=0.9$}\\
				\cline{2-3} \cline{4-5} \cline{6-7}
				&$\gamma=2$ &$\gamma=4$
				&$\gamma=2$ &$\gamma=4$
				&$\gamma=2$ &$\gamma=4$\\
				\midrule
				100	 &8.78   &113.78   &6.40  &62.42   &4.67 &34.24 \\
				200	 &12.42  &212.32   &8.95  &115.50  &6.46 &62.84 \\
				400	 &17.57  &396.20   &12.57 &214.37  &9.00 &116.01 \\
			\end{tabular*}
		\end{spacing}
		{\rule{\temptablewidth}{0.5pt}}
	\end{center}
\end{table}

We are to emphasize that the above procedure
provides a novel technique to verify the positive definiteness
of the discrete convolution kernels,
especially when the first condition
of the criterion \eqref{condition:nonunifrom real quadratic} fails partly.
Recall that $\lambda^{(\mathrm{h})}_{\min}$ denotes
the minimum eigenvalue of the real quadratic form (matrix)
associated with the discrete  L1$_\mathrm{h}$ kernels.
\lan{Theorem \ref{thm:L1h kernel positive definite} is supported by
 Tables \ref{table:L1h kernel verification}-\ref{table:L1h kernel verification on graded meshes}, where the values of $\lambda^{(\mathrm{h})}_{\min}$ are recorded on  the random time meshes with $T=1$ and the graded mesh for three different fractional orders $\alpha=0.1$, $0.5$ and $0.9$. They suggest that Theorem \ref{thm:L1h kernel positive definite} still has a lot of room for improvement, at least on graded meshes.}


\subsection{Analysis of L1$_\mathrm{a}$ kernels}

Due to Theorem \ref{thm:L1 kernel positive definite},
the L1$_\mathrm{a}$ kernels $\La{j}{n}$ in \eqref{def:L1a formula kernel}
would be expected to be positive definite since they are nothing but the averaged version
of L1 kernels $a_{j}^{(n)}$.
Nonetheless, it is invalid.

At first, the first two kernels $\La{0}{n}$ and $\La{1}{n}$ do not maintain the monotonicity property.
According to the definition \eqref{def:L1a formula kernel}, the first two kernels satisfy
\begin{align*}
\La{0}{n}-\La{1}{n}
=\La{0}{n}\kbrab{1-\brat{1+r_n}^{1-\alpha}r_n^{\alpha}-r_n^{\alpha}+r_n}.
\end{align*}
Apparently, $\La{0}{n}<\La{1}{n}$ as the fractional order
$\alpha\rightarrow0$.
In the fractional order limit $\alpha\rightarrow1$, we find $\La{0}{n}>\La{1}{n}$ if the time-step
ratio $r_n<1$, and $\La{0}{n}<\La{1}{n}$ if $r_n>1$.
Always, the integral mean-value theorem gives the following result.

\begin{lemma}
The discrete L1$_\mathrm{a}$ kernels $\La{j}{n}$ in \eqref{def:L1a formula kernel}
satisfy
\begin{align*}
\La{1}{n}>\La{2}{n}>\cdots>\La{n-1}{n}>0\quad\text{but}\quad
\La{0}{n}\ngeqslant \La{1}{n}\quad\text{for $n\ge2$.}
\end{align*}
\end{lemma}

As done in the above subsection, one may remedy this issue by introducing the following
auxiliary kernels
\begin{align}\label{def:auxiliary L1a kernel}
\aLa{a}{0}{n}:=2\La{0}{n}\quad
\text{and}\quad \aLa{a}{j}{n}:=\La{j}{n}
 \quad\text{for $1\le j\le n-1$.}
\end{align}
By the inequality $1+\beta z>(1+z)^\beta$ for $z>0$ and $\beta\in(0,1)$, it is not difficult to check that
\begin{align}\label{inequality:first two L1a kernel}
\aLa{a}{0}{n}-\aLa{a}{1}{n}
&=\La{0}{n}\kbra{2-\brat{1+r_n}^{1-\alpha}r_n^{\alpha}-r_n^{\alpha}+r_n}\nonumber\\
&=\La{0}{n}r_n\kbra{1+2/r_n-\brat{1+1/r_n}^{1-\alpha}-r_n^{\alpha-1}}\nonumber\\
&>\La{0}{n}\brat{1+\alpha-r_n^{\alpha}}.
\end{align}
As seen, a step-ratios restriction $0<r_n\le\sqrt[\alpha]{1+\alpha}$ is necessary
to recover the decreasing property.

To establish the positive definiteness by \cite[Theorem 1.1]{LiaoTangZhou:2020Positive},
we need to confirm that the auxiliary discrete kernels $\aLa{a}{j}{n}$ fulfill
the last two algebraic conditions in \eqref{condition:nonunifrom real quadratic}.
As done in the proof of Lemma \ref{lem:L1h auxiliary  kernel positive definiteness},
one can introduce the following sequence
\begin{align*}
\psi_{n-k}^{\bra{\mathrm{a},n}}
:=\frac{\aLa{a}{n-k}{n}}{\aLa{a}{n-1-k}{n-1}}
\quad \text{for $1\le k\le n-1$}.
\end{align*}
By direct calculations, we have
\begin{align*}
\psi_{1}^{\bra{\mathrm{a},n}}
=\frac12\kbrab{\bra{r_n+1}^{1-\alpha}-r_n^{1-\alpha}+1}
\end{align*}
and
\begin{align*}
\psi_{2}^{\bra{\mathrm{a},n}}=\frac{\bra{r_nr_{n-1}+r_{n-1}+1}^{1-\alpha}
-\bra{r_nr_{n-1}+r_{n-1}}^{1-\alpha}-1}
{\bra{r_{n-1}+1}^{1-\alpha}-r_{n-1}^{1-\alpha}+1}+1.
\end{align*}
Evidently, it is seen that
$\psi_{1}^{\bra{\mathrm{a},n}}\rightarrow1/2$
and $\psi_{2}^{\bra{\mathrm{a},n}}\rightarrow0$
as the fractional order $\alpha\rightarrow1$;
and $\psi_{1}^{\bra{\mathrm{a},n}}\rightarrow1$
and $\psi_{2}^{\bra{\mathrm{a},n}}\rightarrow1$
as $\alpha\rightarrow0$.
So $\psi_{2}^{\bra{\mathrm{a},n}}\ngeqslant\psi_{1}^{\bra{\mathrm{a},n}}$ for $\alpha\in(0,1)$.
Reminding these facts, one can
follow the proof of Lemma \ref{lem:L1h auxiliary  kernel positive definiteness} to prove the following lemma.
It implies that the auxiliary kernels technique
fails to verify the positive definiteness of the L1$_\mathrm{a}$ kernels $\La{j}{n}$,
because the auxiliary kernels $\aLa{a}{j}{n}$ do not fulfill
the third algebraic condition in \eqref{condition:nonunifrom real quadratic}.

\begin{lemma}\label{lem:auxiliary L1a kernel property}
Let $n\ge3$. For the auxiliary L1$_\mathrm{a}$ kernels $\aLa{a}{j}{n}$
in \eqref{def:auxiliary L1a kernel}, it holds that
\begin{align*}
\aLa{a}{j-1}{n-1}\aLa{a}{j+1}{n}\ge\aLa{a}{j}{n-1}\aLa{a}{j}{n}
\quad\text{for $2\le j\le n-2$}\quad\text{but}\quad
\aLa{a}{0}{n-1}\aLa{a}{2}{n}\ngeqslant\aLa{a}{1}{n-1}\aLa{a}{1}{n}.
\end{align*}
\end{lemma}

The above arguments do not negate the positive definiteness
of the L1$_\mathrm{a}$ kernels $\La{j}{n}$ in \eqref{def:L1a formula kernel},
while the numerical computations do. Recall that $\lambda^{(\mathrm{a})}_{\min}$ represents the minimum
eigenvalue of the real quadratic form (matrix) associated with the discrete L1$_\mathrm{a}$ kernels.
\lan{Tables \ref{table:min eigenvalue L1a kernel}-\ref{table:min eigenvalue L1a kernel on graded meshes} record the values of $\lambda^{(\mathrm{a})}_{\min}$ on the time meshes with some fixed step-ratios (more results for other cases of $r_n < 1$ are omitted for brevity) and the graded meshes $t_k=T(k/N)^{\gamma}$ with $T=1$, respectively.}
We observe that the L1$_\mathrm{a}$ kernels are non-positive definite if $r_n>1$, while they may be positive definite if the step-ratios $r_n\le1$.
Up to now, no theoretical proof is available for the later case.

\begin{table}[htb!]
\begin{center}
\caption{The minimum eigenvalue $\lambda_{\min}^{(\mathrm{a})}$ on different time meshes.}
\label{table:min eigenvalue L1a kernel} \vspace*{0.3pt}
\def\temptablewidth{0.85\textwidth}
{\rule{\temptablewidth}{0.5pt}}
\begin{spacing}{1.25}
\begin{tabular*}{\temptablewidth}{@{\extracolsep{\fill}}c|cc|cc|cc}
  \multirow{2}{*}{$N$} &\multicolumn{2}{c|}{$\alpha=0.1$}
  &\multicolumn{2}{c|}{$\alpha=0.5$} &\multicolumn{2}{c}{$\alpha=0.9$}\\
  \cline{2-3} \cline{4-5} \cline{6-7}
         &$r_n=1$  &$r_n=1.1$    &$r_n=1$  &$r_n=1.1$    &$r_n=1$  &$r_n=1.1$\\
  \midrule
  100	 &7.04e-05 &-6.86e-03	 &2.60e-03 &-3.67e+00	 &2.87e-02 &-6.54e+02\\	 	
  200	 &1.90e-05 &-1.78e-02	 &9.27e-04 &-4.30e+02	 &1.35e-02 &-3.48e+06\\	 	
  400	 &5.10e-06 &-4.24e+01	 &3.29e-04 &-5.93e+06	 &6.34e-03 &-9.81e+13\\	
\end{tabular*}
\end{spacing}
{\rule{\temptablewidth}{0.5pt}}
\end{center}
\end{table}

\begin{table}[htb!]
	\begin{center}
		\caption{\lan{The minimum eigenvalue $\lambda_{\min}^{(\mathrm{a})}$ on graded time meshes.}}
		\label{table:min eigenvalue L1a kernel on graded meshes} \vspace*{0.3pt}
		\def\temptablewidth{0.85\textwidth}
		{\rule{\temptablewidth}{0.5pt}}
		\begin{spacing}{1.25}
		\begin{tabular*}{\temptablewidth}{@{\extracolsep{\fill}}c|cc|cc|cc}
			\multirow{2}{*}{$N$} &\multicolumn{2}{c|}{$\alpha=0.1$}
			&\multicolumn{2}{c|}{$\alpha=0.5$} &\multicolumn{2}{c}{$\alpha=0.9$}\\
			\cline{2-3} \cline{4-5} \cline{6-7}
			&$\gamma=2$ &$\gamma=4$
			&$\gamma=2$ &$\gamma=4$
			&$\gamma=2$ &$\gamma=4$\\
			\midrule
			100	 &-1.42e-02	 &-1.65e-01   &-2.98e+00  &-1.06e+03  &-1.96e+02	 &-2.42e+06  \\
			200	 &-1.63e-02  &-2.18e-01	  &-5.97e+00  &-4.22e+03  &-6.81e+02	 &-2.94e+07  \\
			400	 &-1.87e-02	 &-2.88e-01   &-1.19e+01  &-1.69e+04  &-2.37e+03	 &-3.56e+08  \\
		\end{tabular*}
	\end{spacing}
		{\rule{\temptablewidth}{0.5pt}}
	\end{center}
\end{table}

As a special case, we consider
the auxiliary L1$_\mathrm{a}$ kernels $\aLa{a}{j}{n}=\mathsf{a}_j^{(\mathrm{a})}$
on the uniform time mesh. By the definition \eqref{def:L1a formula kernel},
it is not difficult to check that
$$\mathsf{a}_0^{(\mathrm{a})}-2\mathsf{a}_1^{(\mathrm{a})}+\mathsf{a}_2^{(\mathrm{a})}
=a_0^{(\mathrm{a})}-\brab{a_0^{(\mathrm{a})}+a_1^{(\mathrm{a})}}+\frac12\brab{a_1^{(\mathrm{a})}+a_2^{(\mathrm{a})}}
=-\frac12\brab{a_1^{(\mathrm{a})}-a_2^{(\mathrm{a})}}<0.$$
We see that the third condition of
\eqref{condition:unifrom real quadratic} is also not satisfied.
Thus the L\'opez-Marcos criteria in \cite[Proposition 5.2]{Lopez:1990ADifferenceScheme}
 are not enough to ensure the positive definiteness of
the auxiliary L1$_\mathrm{a}$ kernels $\mathsf{a}_j^{(\mathrm{a})}$
and the original L1$_\mathrm{a}$ kernels $a_j^{(\mathrm{a})}$ as well.


\section{Energy dissipation laws of L1-type schemes}
\setcounter{equation}{0}

In this section, the discrete energy stabilities of the proposed
L1-type schemes \eqref{scheme:fully implicit L1 TFCH}-\eqref{scheme:fully implicit L1h TFCH}
are established by making use of the above theoretical results on the L1 and L1$_\mathrm{h}$ kernels.
Here and hereafter,
we use the standard norms of the Sobolev space $H^{m}\bra{\Omega}$
and the $L^p\bra{\Omega}$ space.
For any functions $v$ and $w$ belonging to the zero-mean space
$\mathbb{\mathring V}:=\big\{v\in L^2\bra{\Omega}\,|\, \brab{v,1}=0\big\}$,
the $H^{-1}$ inner product $\brab{v,w}_{-1}:=\brab{\brat{-\Delta}^{-1}v,w}$ and
the induced norm $\mynormb{v}_{-1}:=\sqrt{\brab{v,v}_{-1}}$
may be used, which will not be introduced specifically.
\subsection{Variable-step L1 scheme}
At first, we investigate the discrete volume conservation property
and unique solvability of the variable-step L1
scheme \eqref{scheme:fully implicit L1 TFCH}.

\begin{lemma}\label{lem:TFCH L1 volume law}
The variable-step L1 scheme \eqref{scheme:fully implicit L1 TFCH} conserves the volume,
\begin{align*}
\brab{\phi^n,1}=\brab{\phi^{n-1},1}\quad \text{for $1\le n\le N$.}
\end{align*}
\end{lemma}

\begin{proof}
Taking the inner product of \eqref{scheme:fully implicit L1 TFCH} with 1, one applies
 the Green's formula to find
\begin{align*}
\brab{\bra{\partial_{\tau}^{\alpha}\phi}^n,1}
=\kappa\brab{\Delta \mu^n,1}=0.
\end{align*}
Multiplying both sides of the equality by $\theta_{m-n}^{(m)}$ and summing up from $n=1$ to $m$, we have
\begin{align*}
\braB{\sum_{n=1}^m\theta_{m-n}^{(m)}\bra{\partial_{\tau}^{\alpha}\phi}^n,1}
=\braB{\sum_{n=1}^m\theta_{m-n}^{(m)}\sum_{j=1}^na_{n-j}^{(n)}\diff\phi^j,1}=0\quad\text{for $m\ge1$.}
\end{align*}
By exchanging the summation order and applying the discrete orthogonal identity \eqref{equality:orthogonal identity},
it arrives at $\brab{\diff\phi^m,1}=0$ for $m\ge1$.
The assertion follows and the proof is completed.
\end{proof}

\begin{theorem}\label{thm:L1 convexity solvability TFCH}
Under the time-step restriction
\begin{align}\label{inequality:L1 time step condition solvability TFCH}
\tau_n\le \sqrt[\alpha]{\frac{4\epsilon^2}{\kappa\Gamma(2-\alpha)}},
\end{align}
the variable-step L1 scheme \eqref{scheme:fully implicit L1 TFCH} is uniquely solvable.
\end{theorem}
\begin{proof}
For any fixed time-level indexes $n\ge1$,
we consider the following energy functional $G[z]$ on the space
$\mathbb{V}_{h}^{*}:=\big\{z\in\mathbb{V}_{h}\,|\, \brab{z,1}=\brab{\phi^{n-1},1}\big\},$
\begin{align*}
G[z]:=\frac{a_0^{(n)}}{2}\mynormb{z-\phi^{n-1}}_{-1}^2
+\brab{\mathcal{L}^{n-1},z-\phi^{n-1}}_{-1}
+\frac{\epsilon^2}{2}\kappa\mynormb{\nabla z}^2
+\frac\kappa4\mynormb{z}_{L^4}^4-\frac{\kappa}{2}\mynormb{z}^2,
\end{align*}
where we use the notation
$\mathcal{L}^{n-1}:=\sum_{k=1}^{n-1}a_{n-k}^{(n)}\diff \phi^{k}$
for brevity.
The time-step restriction \eqref{inequality:L1 time step condition solvability TFCH}
implies that the discrete L1 kernel $a_0^{(n)}\ge\kappa/\brat{4\epsilon^2}$.
By using the inequality
$$\mynormb{v}^2\le \mynormb{\nabla v}\mynormb{v}_{-1}
\le \epsilon^2\mynormb{\nabla v}^2
+\frac{1}{4\epsilon^2}\mynormb{v}_{-1}^2\quad
\text{for any $v\in\mathbb{\mathring V}$},$$
we see that the energy functional $G[z]$
is convex with respect to $z$, that is,
\begin{align*}
\frac{\zd^2G}{\zd s^2}[z+s\psi]\Big|_{s=0}
&=a_0^{(n)}\mynormb{\psi}_{-1}^2
+\kappa\epsilon^2\mynormb{\nabla \psi}^2
-\kappa\mynormb{\psi}^2
+3\kappa\mynormb{z\psi}^2\\
&\ge \brab{a_0^{(n)}-\frac{\kappa}{4\epsilon^2}}\mynormb{\psi}_{-1}^2
+3\kappa\mynormb{z\psi}^2>0.
\end{align*}
It is easily to show
that the functional $G[z]$ is coercive on $\mathbb{V} ^{*}$, that is,
\begin{align*}
G[z]&\ge \frac{\kappa}{4}\mynorm{z}_{L^4}^4-\frac{\kappa}{2}\mynormb{z}^2
-\frac{1}{2a_0^{(n)}}\mynormb{\mathcal{L}^{n-1}}^2
\ge\frac{\kappa}{2}\mynormb{z}^2
-\frac{1}{2a_0^{(n)}}\mynormb{\mathcal{L}^{n-1}}^2
-\kappa\abs{\Omega },
\end{align*}
where the inequality
$\mynorm{v}_{L^4}^4\ge 4\mynormb{v}^2-4\abs{\Omega }$
has been used in the last step.
So the functional $G[z]$ has a unique minimizer,
which implies the L1 scheme \eqref{scheme:fully implicit L1 TFCH} exists a unique solution.
\end{proof}

\begin{remark}
Let the fractional order $\alpha\rightarrow 1$, the variable-step
L1 scheme \eqref{scheme:fully implicit L1 TFCH} approaches
 the standard backward Euler scheme
\begin{align}\label{scheme:backward Euler CH}
\partial_\tau\phi^n=\kappa\Delta \mu^n\quad\text{with}\quad
\mu^n=\bra{\phi^n}^3-\phi^n-\epsilon^2\Delta \phi^n,
\quad  n\ge 1,
\end{align}
which is uniquely solvable under the time-step restriction
$\tau_n\le 4\epsilon^2/\kappa$, see
\cite[Theorem 2.2]{Xu:2019OnFullImplicitScheme} with our notation $\kappa=1/\epsilon$.
The time-step condition  \eqref{inequality:L1 time step condition solvability TFCH}
is asymptotically compatible with the above restriction in the fractional order limit $\alpha\rightarrow 1$.
\end{remark}

Let $E\kbra{\phi^{n}}$ be the discrete version of the
free energy functional \eqref{cont:classical free energy},
\begin{align}\label{def:discrete original energy}
E\kbra{\phi^{n}}:=\frac{\epsilon^2}{2}\mynormb{\nabla \phi^n}^2
+\brab{F(\phi^n),1}\quad\text{with}\quad
F(\phi^n):=\frac14\brab{\brat{\phi^n}^2-1}^2\quad\text{for $n\ge 0$.}
\end{align}
The discrete counterpart $\mathcal{E}_\alpha\kbra{\phi^{n}}$
of the variational energy functional \eqref{cont:TFCH variational energy} is given by
\begin{align*}
\mathcal{E}_\alpha\kbrab{\phi^0}:=E\kbrab{\phi^0}\quad \text{and}\quad
\mathcal{E}_\alpha\kbrab{\phi^{n}}:=E\kbrab{\phi^{n}}
+\frac{\kappa}{2}\sum_{j=1}^np_{n-j}^{(n)}\mynormb{\nabla \mu^j}^2
\quad\text{for $n\ge 1$,}
\end{align*}
where the DCC kernels $p_{n-j}^{(n)}$ with respect to the L1 kernels
$a_{n-j}^{(n)}$ are used to simulate the Riemann-Liouville fractional integral
$\bra{\mathcal{I}_t^\alpha v}(t_n)\approx\sum_{j=1}^np_{n-j}^{(n)}v^j$,
cf. \cite{LiaoTangZhou:2020EnergyMaximumBound,LiaoYanZhang:2018Unconditional}.
\begin{theorem}\label{thm:L1 energy law TFCH}
Under the time-step restriction \eqref{inequality:L1 time step condition solvability TFCH},
the variable-step L1 scheme \eqref{scheme:fully implicit L1 TFCH}
preserves the following discrete energy dissipation law
\begin{align*}
\partial_\tau\mathcal{E}_\alpha\kbra{\phi^n}\le 0 \quad \text{for $1\le n\le N$.}
\end{align*}
\end{theorem}

\begin{proof}
Making the inner product of the equation \eqref{scheme:fully implicit L1 TFCH} by  $\bra{-\Delta }^{-1}\diff \phi^{n}/\kappa$, one obtains
\begin{align}\label{thmProof:L1 energy law inner}
\frac{1}{\kappa}\brab{\bra{\partial_\tau^\alpha \phi}^n,\diff \phi^n}_{-1}
+\brab{\bra{\phi^n}^3-\phi^{n},\diff \phi^n}
-\brab{\epsilon^2\Delta \phi^n,\diff \phi^n}=0.
\end{align}
An application of the inequality
$$\brat{a^3-a}\bra{a-b}\ge \frac14\bra{a^2-1}^2-\frac14\bra{b^2-1}^2-\frac12\bra{a-b}^2$$
to the second term of equation \eqref{thmProof:L1 energy law inner} yields
\[
\brab{f\bra{\phi^n},\diff \phi^n}
\ge \brab{F(\phi^n),1}-\brab{F(\phi^{n-1}),1}
-\frac12\mynormb{\diff \phi^n}^2.
\]
For the third term of equation \eqref{thmProof:L1 energy law inner},
the identity $2a(a-b)=a^2-b^2+(a-b)^2$ gives
\begin{align*}
-\brab{\epsilon^2\Delta \phi^n,\diff \phi^n}
=\frac{\epsilon^2}{2}\mynormb{\nabla \phi^n}^2
-\frac{\epsilon^2}{2}\mynormb{\nabla \phi^{n-1}}^2
+\frac{\epsilon^2}{2}\mynormb{\diff\nabla\phi^n}^2.
\end{align*}
Substituting the above results into equation
\eqref{thmProof:L1 energy law inner}, one has
\begin{align}\label{thmProof:L1 energy law norm}
\frac{1}{\kappa}\brab{\bra{\partial_\tau^\alpha \phi}^n,\diff \phi^n}_{-1}
+\frac{\epsilon^2}{2}\mynormb{\diff \nabla \phi^n}^2
-\frac12\mynormb{\diff \phi^n}^2
+E\kbra{\phi^n}
\le E\kbra{\phi^{n-1}}.
\end{align}
For the first term  of  \eqref{thmProof:L1 energy law norm},
the first inequality in Theorem \ref{thm:L1 kernel positive definite} yields
\begin{align*}
\frac{1}{\kappa}\brab{\bra{\partial_\tau^\alpha \phi}^n,\diff \phi^n}_{-1}
\ge \frac{\kappa}{2}\sum_{j=1}^np_{n-j}^{(n)}\mynormb{\nabla \mu^j}^2
-\frac{\kappa}{2}\sum_{j=1}^{n-1}p_{n-1-j}^{(n-1)}\mynormb{\nabla \mu^j}^2
+\frac{a_0^{(n)}}{2\kappa}\mynormb{\diff\phi^n}_{-1}^2,
\end{align*}
where the following identity has been used in the above derivation
\begin{align*}
\mynormB{\sum_{\ell=1}^ja_{j-\ell}^{(j)}\diff\phi^\ell}_{-1}^2
=\braB{\sum_{\ell=1}^ja_{j-\ell}^{(j)}\diff\phi^\ell,
\bra{-\Delta }^{-1}\sum_{\ell=1}^ja_{j-\ell}^{(j)}\diff\phi^\ell}
=\mynormb{\kappa\nabla \mu^j}^2.
\end{align*}
Furthermore, we have
$$-\frac12\mynormb{\diff\phi^n}^2\ge -\frac{1}{8\epsilon^2}\mynormb{\diff\phi^n}_{-1}^2
-\frac{\epsilon^2}{2}\mynormb{\diff\nabla \phi^n}^2.$$
Inserting the above estimates into the left hand side of \eqref{thmProof:L1 energy law norm},
one gets
\begin{align*}
\frac1{2\kappa}\braB{a_0^{(n)}-\frac{\kappa}{4\epsilon^2}}\mynormb{\diff\phi^n}_{-1}^2
+\mathcal{E}_\alpha\kbrab{\phi^n}
\le \mathcal{E}_\alpha\kbrab{\phi^{n-1}}.
\end{align*}
Then the claimed result follows from the time-step condition
\eqref{inequality:L1 time step condition solvability TFCH} immediately.
\end{proof}

\begin{remark}\label{remark:L1 energy compatibility}
Under the restriction $\tau_n\le 4\epsilon^2/\kappa$,
the backward Euler scheme \eqref{scheme:backward Euler CH} for the classical CH model
preserves the energy dissipation law
\cite[Theorem 2.2]{Xu:2019OnFullImplicitScheme},
\begin{align}\label{inequality:backward Euler energy decay law}
\partial_\tau E\kbra{\phi^n}
+\frac{\kappa}{2}\mynormb{\nabla \mu^n}^2\le 0
\quad \text{for  $1\le n\le N$.}
\end{align}
As the fractional index $\alpha\rightarrow 1$,
the definition \eqref{scheme:L1 formula} shows that
the L1 kernels
$a_0^{(n)}\rightarrow 1/\tau_n$ and $a_{n-k}^{(n)}\rightarrow 0$ for $1\le k\le n-1$.
Corresponding, the DOC kernels $\theta_0^{(n)}\rightarrow \tau_n$
and $\theta_{n-k}^{(n)}\rightarrow 0$ for $1\le k\le n-1$,
and the DCC kernels $p_{n-k}^{(n)}\rightarrow \tau_k$ for $1\le k\le n$.
So the variational energy dissipation law in Theorem \ref{thm:L1 energy law TFCH}
is asymptotically compatible with \eqref{inequality:backward Euler energy decay law}
in the sense that
\[
\partial_\tau\mathcal{E}_\alpha\kbra{\phi^n}\le 0\quad\longrightarrow\quad
\partial_\tau E\kbra{\phi^n}
+\frac{\kappa}{2}\mynormb{\nabla \mu^n}^2\le 0
\quad\text{as $\alpha\rightarrow 1$.}
\]
\end{remark}

Theorem \ref{thm:L1 energy law TFCH} implies that the solution of the
L1 scheme \eqref{scheme:fully implicit L1 TFCH}
is bounded in the $H^1$ norm.

\begin{corollary}\label{corollary:L1 H1 solution bound}
The solution of the variable-step L1 scheme
\eqref{scheme:fully implicit L1 TFCH} satisfies,
\[
\mynormb{\phi^n}_{H ^1}
\le\sqrt{\bra{4E\kbra{\phi^0}+\bra{2\epsilon^2
+\epsilon^4}\absb{\Omega }}/(2\epsilon^2)}:=c_0
\quad\text{for $n\ge 1$,}
\]
where the constant $c_0$ is dependent on the domain $\Omega$,
the parameter $\epsilon$ and the initial value $\phi^0$,
but independent of the time $t_n$, step sizes $\tau_n$ and time-step ratios $r_n$.
\end{corollary}
\begin{proof}
The discrete energy law in Theorem \ref{thm:L1 energy law TFCH} gives
$E[\phi^{0}]\ge \mathcal{E}_\alpha[\phi^n]\ge E[\phi^n]$.
Then by the following inequality
$$\mynorm{v}_{L^4}^4\ge (2+2\epsilon^2)\mynorm{v}^2-(1+\epsilon^2)^2\absb{\Omega},$$
one has
\begin{align*}
4E\kbra{\phi^0}&\ge 2\epsilon^2\mynormb{\nabla \phi^n}^2
+\mynorm{\phi^n}_{L^4}^4-2\mynormb{\phi^n}^2+\absb{\Omega}\\
&\ge 2\epsilon^2\mynormb{\nabla \phi^n}^2
+2\epsilon^2\mynormb{\phi^n}^2-\bra{2\epsilon^2+\epsilon^4}\absb{\Omega },
\end{align*}
which yields the claimed solution bound immediately.
This completes the proof.
\end{proof}

By using this type solution bound, an $L^2$ norm error estimate
for the L1 scheme \eqref{scheme:fully implicit L1 TFCH} can be derived
by following the analysis in \cite{LiaoZhuWang:2020TFAC}, but we omit it here for brevity.

\subsection{Variable-step L1$_\mathrm{h}$ scheme}

Now we investigate the volume-conserving property, the unique solvability
and the discrete energy stability for the variable-step L1$_\mathrm{h}$ scheme
\eqref{scheme:fully implicit L1h TFCH}.
By following the proof of Lemma \ref{lem:TFCH L1 volume law} with the corresponding
DOC kernels $\theta_{n-j}^{(\mathrm{h},n)}$ with respect to the original L1$_\mathrm{h}$ kernels $a_{n-j}^{(\mathrm{h},n)}$,
it is easy to obtain the following result.

\begin{lemma}\label{lem:TFCH L1h volume law}
The variable-step L1$_\mathrm{h}$ scheme
\eqref{scheme:fully implicit L1h TFCH} conserves the volume,
\begin{align*}
\brab{\phi^n,1}=\brab{\phi^{n-1},1}\quad \text{for $1\le n\le N$.}
\end{align*}
\end{lemma}
Consider a discrete energy functional $G_{\mathrm{h}}[z]$ defined on the
volume-conserving space $\mathbb{V}_{h}^{*}$ as
\begin{align*}
G_{\mathrm{h}}[z]&:=\frac{1}{2}\Lh{0}{n}\mynormb{z-\phi^{n-1}}_{-1}^2
+\brab{\mathcal{L}_{\mathrm{h}}^{n-1},z-\phi^{n-1}}_{-1}
+\frac{\epsilon^2}{4}\kappa\mynormb{\nabla \brab{z+\phi^{n-1}}}^2\\
&\quad+\frac\kappa8\mynormb{z}_{L^4}^4
+\frac\kappa4\brab{\brab{\phi^{n-1}}^2,z^2}
+\frac\kappa6\brab{\brab{\phi^{n-1}}^3,z}
-\frac\kappa4\mynormb{\brab{z+\phi^{n-1}}}^2,
\end{align*}
where the notation
$\mathcal{L}_{\mathrm{h}}^{n-1}:=\sum_{k=1}^{n-1}\Lh{n-k}{n}\diff \phi^{k}$.
By following the convexity argument performed in the proof of
Lemma \ref{thm:L1 convexity solvability TFCH},
it is not difficult to prove the unique solvability of
\eqref{scheme:fully implicit L1h TFCH}.

\begin{theorem}
Under the time-step restriction
\begin{align}\label{inequality:L1h solvability TFCH}
\tau_n\le 2\sqrt[\alpha]{\frac{4\epsilon^2}{\kappa\Gamma(2-\alpha)}},
\end{align}
the variable-step L1$_{\mathrm{h}}$ scheme
\eqref{scheme:fully implicit L1h TFCH} is uniquely solvable.
\end{theorem}

\begin{remark}
Consider the following  Crank-Nicolson scheme for the CH model 
\begin{align}\label{scheme:CN type classical CH}
\partial_\tau\phi^n=\kappa\Delta \mu^{n-\frac12}\quad\text{with}\quad
\mu^{n-\frac12}=f\brat{\phi}^{n-\frac12}-\epsilon^2\Delta \phi^{n-\frac12},
\quad  n\ge 1.
\end{align}
It is not difficult to check that it is uniquely solvable under the time-step restriction
$\tau_n\le 8\epsilon^2/\kappa$. As the fractional order $\alpha\rightarrow 1$,
the definition \eqref{scheme:L1h formula} of
the original L1$_{\mathrm{h}}$ kernels $\Lh{n-k}{n}$ implies that
$$\Lh{0}{n}\rightarrow 1/\tau_n\quad
\text{and}\quad \Lh{n-k}{n}\rightarrow 0 \quad\text{for $1\le k\le n-1$.}$$
Thus the variable-step L1$_{\mathrm{h}}$ scheme \eqref{scheme:fully implicit L1h TFCH} degenerates into
the Crank-Nicolson scheme \eqref{scheme:CN type classical CH}.
We see that the time-step condition \eqref{inequality:L1h solvability TFCH}
for the variable-step L1$_\mathrm{h}$ scheme
\eqref{scheme:fully implicit L1h TFCH} is sharp in the sense that
it approaches the time-step restriction for \eqref{scheme:CN type classical CH} in the limit $\alpha\rightarrow 1.$
\end{remark}

By virtues of Theorem \ref{thm:L1h kernel positive definite}
for the original discrete kernels $\Lh{n-k}{n}$,
we are to build up a discrete variational energy dissipation law for
the L1$_\mathrm{h}$ scheme \eqref{scheme:fully implicit L1h TFCH}.
As the main difference to the above case for the L1 scheme, Theorem \ref{thm:L1h kernel positive definite}
involves the auxiliary L1$_\mathrm{h}$ kernels $\aLh{a}{n-j}{n}$ and the associated DCC kernels $\aLh{p}{n-j}{n}$.
We define the following (unusual) discrete variational energy $\mathcal{E}_\alpha^{(\mathrm{h})}$
\begin{align*}
\mathcal{E}_\alpha^{(\mathrm{h})}\kbra{\phi^0}:=E\kbra{\phi^0}\quad \text{and}\quad
\mathcal{E}_\alpha^{(\mathrm{h})}\kbra{\phi^{n}}:=E\kbra{\phi^{n}}
+\frac{1}{2\kappa}\sum_{j=1}^n\aLh{p}{n-j}{n}\mynormB{\sum_{\ell=1}^j\aLh{a}{j-\ell}{j}\diff\phi^\ell}_{-1}^2
\quad\text{for $n\ge 1$,}
\end{align*}
where the original energy $E\kbra{\phi^{n}}$
is defined in \eqref{def:discrete original energy}.

\begin{theorem}\label{thm:L1h energy law TFCH}
The variable-step L1$_\mathrm{h}$ scheme \eqref{scheme:fully implicit L1h TFCH}
is unconditionally energy stable in the sense that it
preserves the following discrete energy dissipation law
\begin{align*}
\partial_{\tau}\mathcal{E}_\alpha^{(\mathrm{h})}\kbra{\phi^n}
\le 0\quad \text{for $1\le n\le N$.}
\end{align*}
\end{theorem}
\begin{proof}
Taking the inner product of the equation \eqref{scheme:fully implicit L1h TFCH} by  $\bra{-\Delta }^{-1}\diff \phi^{n}/\kappa$, one gets
\begin{align}\label{thmProof:L1h energy law inner}
\frac{1}{\kappa}\braB{\bra{\partial_{\mathrm{h}\tau}^\alpha\phi}^{n-\frac12},
\diff \phi^n}_{-1}
+\brab{f\brat{\phi}^{n-\frac12},\diff \phi^n}
-\brab{\epsilon^2\Delta\phi^{n-\frac12},\diff \phi^n}=0.
\end{align}
For the first term,
the first inequality in Theorem \ref{thm:L1h kernel positive definite} gives
\begin{align*}
\braB{\bra{\partial_{\mathrm{h}\tau}^\alpha\phi}^{n-\frac12},
\diff \phi^n}_{-1}
\ge \frac{1}{2}\sum_{j=1}^n\aLh{p}{n-j}{n}
\mynormB{\sum_{\ell=1}^j\aLh{a}{j-\ell}{j}\diff\phi^\ell}_{-1}^2
-\frac{1}{2}\sum_{j=1}^{n-1}\aLh{p}{n-1-j}{n-1}
\mynormB{\sum_{\ell=1}^j\aLh{a}{j-\ell}{j}\diff\phi^\ell}_{-1}^2.
\end{align*}
For the second term of \eqref{thmProof:L1h energy law inner},
it follows from \cite[Appendix A]{LiaoTangZhou:2020EnergyMaximumBound}
that
\[
\brab{f\brat{\phi}^{n-\frac12},\diff \phi^n}
=\brab{F(\phi^n),1}-\brab{F(\phi^{n-1}),1}
+\frac{1}{12}\mynormb{\diff\phi^n}_{L^4}^4.
\]
For the third term of \eqref{thmProof:L1h energy law inner}, one has
\[
-\brab{\epsilon^2\Delta\phi^{n-\frac12},\diff \phi^n}
=\frac{\epsilon^2}{2}\mynormb{\nabla\phi^n}^2
-\frac{\epsilon^2}{2}\mynormb{\nabla\phi^{n-1}}^2.
\]
Inserting the above results into the equation \eqref{thmProof:L1h energy law inner}
yields the discrete energy dissipation law immediately.
This completes the proof.
\end{proof}

\begin{remark}\label{remark:L1h energy compatibility}
It is not difficulty to derive that the Crank-Nicolson scheme
\eqref{scheme:CN type classical CH} unconditionally
preserves the following discrete energy law, that is,
\[
\partial_\tau E\kbra{\phi^n}
+\kappa\mynormb{\nabla \mu^{n-\frac12}}^2\le 0
\quad \text{for  $1\le n\le N$.}
\]
As the fractional order $\alpha\rightarrow 1$,
the definitions \eqref{scheme:L1h formula} and \eqref{def:auxiliary L1h kernel} give
$\aLh{a}{0}{n}\rightarrow 2/\tau_n$
and $\aLh{a}{n-k}{n}\rightarrow 0$ for $1\le k\le n-1$.
In turn,
the corresponding DOC kernels $\aLh{\uptheta}{0}{n}\rightarrow \tau_n/2$
and $\aLh{\uptheta}{n-k}{n}\rightarrow 0$ for $1\le k\le n-1$,
and the DCC kernels $\aLh{p}{n-k}{n}\rightarrow \tau_k/2$ for $1\le k\le n$.
We have
\begin{align*}
\mathcal{E}_\alpha^{(\mathrm{h})}\kbra{\phi^{n}}
\quad\longrightarrow\quad
E\kbra{\phi^{n}}
+\frac{1}{\kappa}\sum_{j=1}^n\tau_j\mynormb{\partial_{\tau}\phi^j}_{-1}^2\quad\text{as $\alpha\rightarrow 1$.}
\end{align*}
The equation \eqref{scheme:CN type classical CH} gives
$\mynormb{\partial_{\tau}\phi^n}_{-1}=\kappa\mynormb{\nabla \mu^{n-\frac12}}$.
Then it holds that
\[
\partial_\tau\mathcal{E}_\alpha^{(\mathrm{h})}\kbra{\phi^n}\le0\quad\longrightarrow\quad
\partial_\tau E\kbra{\phi^n}
+\kappa\mynormb{\nabla \mu^{n-\frac12}}^2\le 0
\quad\text{as $\alpha\rightarrow 1$,}
\]
which is just the discrete energy dissipation law
of the scheme \eqref{scheme:CN type classical CH} for the CH model.
In this sense, we say that the variational energy dissipation law
in Theorem \ref{thm:L1h energy law TFCH}
is asymptotically compatible in the fractional order limit $\alpha\rightarrow 1.$
\end{remark}
By following the similar fashion in Corollary \ref{corollary:L1 H1 solution bound},
one can derive the following priori estimate
for the variable-step L1$_\mathrm{h}$ scheme
\eqref{scheme:fully implicit L1h TFCH}.
The involved constant $c_0$ is defined in Corollary \ref{corollary:L1 H1 solution bound}.

\begin{corollary}
The solution of the variable-step L1$_\mathrm{h}$ scheme
\eqref{scheme:fully implicit L1h TFCH} can be bounded by
\[
\mynormb{\phi^n}_{H ^1}\le c_0
\quad\text{for $n\ge 1$.}
\]
\end{corollary}

\section{Numerical experiments}
In this section, we examine the performance
of the variable-step methods
\eqref{scheme:fully implicit L1 TFCH}-\eqref{scheme:fully implicit L1a TFCH}
for the TFCH equation.
The Fourier pseudo-spectral method is employed for the spatial discretization
\cite{ChengWangWise:2016WeaklyEnergy,ChengWangWise:2019SPFC}.
The resulting nonlinear system at each time level
is solved by using a simple fixed-point iteration with the termination error $10^{-12}$.
The sum-of-exponentials technique \cite{Jiang:2017SOE}
with an absolute tolerance error $\epsilon=10^{-12}$ and cut-off time $\Delta{t}=\tau_1$
is always adopted in our numerical simulations to reduce the computational cost and storage.

\subsection{Accuracy verification}
\begin{example}
To verify the temporal accuracy,
we solve the TFCH model \eqref{cont:TFCH model} by adding a forcing term
$\partial_{t}^{\alpha}\Phi=\kappa\Delta\mu+g(\mathbf{x},t)$
with the model parameters $\kappa=1$ and $\epsilon=0.5$
for $\mathbf{x}\in(0,2\pi)^{2}$ and $0<t\le1$
such that $\Phi=\omega_{1+\sigma}(t)\sin(x)\sin(y)$
with a regularity parameter $\sigma\in(0,1)$.
\end{example}

Let the final time $T=1$.
We take the graded time mesh $t_k=(k/N_0)^\gamma$ for $0\le k\le N_0$
in the interval $[0, T_0]$, where
$T_0=\min\{1/\gamma,T\}$ and $N_0=\lceil \frac{N}{T+1-\gamma^{-1}}\rceil$.
In the remainder interval $[T_{0},T]$,
the random time meshes
$\tau_{N_{0}+k}:=(T-T_{0})s_{k}/S_1$ for $1\le k\le N_1$
are used
by setting $N_1:=N-N_0$ and $S_1=\sum_{k=1}^{N_1}s_k$,
where $s_k\in(0,1)$ are random numbers.
The spatial domain $\Omega=(0,2\pi)^2$ is discretized by using
$128^2$ uniform grids.
The $L^2$ norm error $e(N):=\max_{1\le{n}\le{N}}\mynorm{\Phi^n-\phi^n}$
is recorded in each run and the experimental order is evaluated by
$$\text{Order}\approx\frac{\log\bra{e(N)/e(2N)}}{\log\bra{\tau(N)/\tau(2N)}},$$
where $\tau(N)$ denotes the maximum time-step size for total $N$ subintervals.
The accuracy tests are performed by taking
the fractional order $\alpha=0.4$, the regularity parameter $\sigma=0.4$
for three grading parameters  $\gamma=3,4$ and 5.
The previous analysis
\cite{LiaoYanZhang:2018Unconditional,JiLiaoZhang:2020Simple,LiaoZhuWang:2020TFAC}
for the L1 formula suggest an optimal graded parameter $\gamma_{\mathrm{opt}}:=(2-\alpha)/\sigma=4$
to achieve the optimal accuracy $O\bra{\tau^{2-\alpha}}$.

\begin{table}[htb!]
\begin{center}
\caption{Numerical accuracy of the  L1 scheme
\eqref{scheme:fully implicit L1 TFCH} for $\alpha=0.4,\,\sigma=0.4$.}
\label{table:L1 error alph 04} \vspace*{0.3pt}
\def\temptablewidth{1.0\textwidth}
{\rule{\temptablewidth}{0.5pt}}
\begin{tabular*}{\temptablewidth}{@{\extracolsep{\fill}}cccccccccc}
\multirow{2}{*}{$N$} &\multirow{2}{*}{$r_{\max}$} &\multicolumn{2}{c}{$\gamma=3$} &\multirow{2}{*}{$r_{\max}$} &\multicolumn{2}{c}{$\gamma=4$} &\multirow{2}{*}{$r_{\max}$}&\multicolumn{2}{c}{$\gamma=5$} \\
             \cline{3-4}          \cline{6-7}         \cline{9-10}
         &          &$e(N)$   &Order&         &$e(N)$    &Order &         &$e(N)$   &Order\\
\midrule
  40     &11.21	 &5.04e-02 &$-$  &15.00	 &1.35e-02  &$-$   &31.00	 &7.12e-03 &$-$\\
  80     &33.05	 &2.19e-02 &1.52 &28.30	 &4.45e-03	&1.63  &36.86	 &2.63e-03	 &1.55\\
  160    &48.79	 &9.54e-03 &1.11 &91.41	 &1.47e-03	&1.80  &448.21	 &9.03e-04	 &1.58\\
  320    &430.56 &4.15e-03 &1.36 &32.54	 &4.88e-04	&1.48  &155.60	 &2.99e-04	 &1.64\\
\end{tabular*}
{\rule{\temptablewidth}{0.5pt}}
\end{center}
\end{table}	

\begin{table}[htb!]
\begin{center}
\caption{Numerical accuracy of the  L1$_{\mathrm{h}}$ scheme
\eqref{scheme:fully implicit L1h TFCH} for $\alpha=0.4,\,\sigma=0.4$.}
\label{table:L1h error alph 04} \vspace*{0.3pt}
\def\temptablewidth{1.0\textwidth}
{\rule{\temptablewidth}{0.5pt}}
\begin{tabular*}{\temptablewidth}{@{\extracolsep{\fill}}cccccccccc}
\multirow{2}{*}{$N$} &\multirow{2}{*}{$r_{\max}$} &\multicolumn{2}{c}{$\gamma=3$} &\multirow{2}{*}{$r_{\max}$} &\multicolumn{2}{c}{$\gamma=4$} &\multirow{2}{*}{$r_{\max}$}&\multicolumn{2}{c}{$\gamma=5$} \\
             \cline{3-4}          \cline{6-7}         \cline{9-10}
         &       &$e(N)$     &Order&         &$e(N)$    &Order &         &$e(N)$   &Order\\
\midrule
  40     &75.84	 &1.02e-02   &$-$  &18.06	 &4.41e-03  &$-$   &31.00	 &9.76e-03 &$-$\\
  80     &29.77	 &4.42e-03	 &1.41 &22.24	 &1.85e-03	&1.33  &107.09	 &3.84e-03 &1.55\\
  160    &23.73	 &1.92e-03	 &1.19 &15.65	 &5.57e-04	&1.52  &151.87	 &1.10e-03 &1.51\\
  320    &79.85	 &8.37e-04	 &1.12 &200.41	 &1.89e-04	&1.60  &39.06	 &3.25e-04 &1.82\\
\end{tabular*}
{\rule{\temptablewidth}{0.5pt}}
\end{center}
\end{table}	

\begin{table}[htb!]
\begin{center}
\caption{Numerical accuracy of the  L1$_{\mathrm{a}}$ scheme
\eqref{scheme:fully implicit L1a TFCH} for $\alpha=0.4,\,\sigma=0.4$.}
\label{table:L1a error alph 04} \vspace*{0.3pt}
\def\temptablewidth{1.0\textwidth}
{\rule{\temptablewidth}{0.5pt}}
\begin{tabular*}{\temptablewidth}{@{\extracolsep{\fill}}cccccccccc}
\multirow{2}{*}{$N$} &\multirow{2}{*}{$r_{\max}$} &\multicolumn{2}{c}{$\gamma=3$} &\multirow{2}{*}{$r_{\max}$} &\multicolumn{2}{c}{$\gamma=4$} &\multirow{2}{*}{$r_{\max}$}&\multicolumn{2}{c}{$\gamma=5$} \\
             \cline{3-4}          \cline{6-7}         \cline{9-10}
         &          &$e(N)$   &Order&        &$e(N)$    &Order &         &$e(N)$   &Order\\
\midrule
  40     &44.98	 &1.21e-02   &$-$  &167.41	 &9.55e-03  &$-$   &31.00	 &9.06e-03 &$-$\\
  80     &14.44	 &5.28e-03	 &1.13 &15.00	 &3.19e-03	&1.69  &104.61	 &4.07e-03 &1.47\\
  160    &42.06	 &2.30e-03	 &1.37 &145.46	 &1.05e-03	&1.77  &31.00	 &1.30e-03 &1.33\\
  320    &86.02	 &1.00e-03	 &1.15 &264.04	 &3.65e-04	&1.32  &48.28	 &4.81e-04 &1.66\\
\end{tabular*}
{\rule{\temptablewidth}{0.5pt}}
\end{center}
\end{table}	

The numerical errors are tabulated in Tables \ref{table:L1 error alph 04}-\ref{table:L1a error alph 04}.
We observe that the accuracies of L1-type schemes
\eqref{scheme:fully implicit L1 TFCH}-\eqref{scheme:fully implicit L1a TFCH}
only reach $O\brab{\tau^{\gamma\sigma}}$
when the graded parameter $\gamma<\gamma_{\text{opt}}$;
while the optimal accuracy $O\brab{\tau^{2-\alpha}}$
can be achieved when the graded parameter
$\gamma\ge\gamma_{\text{opt}}$.
Also, the maximum time-step ratios (denoted by $r_{\max}$)
recorded in Tables \ref{table:L1 error alph 04}-\ref{table:L1a error alph 04}
indicate that the proposed L1-type methods
are robust with respect to the step-size variations.

\subsection{Simulation of coarsening dynamics}

\begin{example}\label{example:coarsen dynamic}
We next simulate the coarsening dynamics of the TFCH equation \eqref{cont:TFCH model}.
The initial condition is taken as
$\phi_0(\mathbf{x})=\text{rand}(\mathbf{x})$,
where $\text{rand}(\mathbf{x})$ generates uniform random numbers
between $-0.001$ to $0.001$.
The mobility coefficient $\kappa=0.01$ and
the interfacial thickness $\epsilon=0.05$.
The spatial domain $(0,2\pi)^2$ is discretized by using $128^2$ spatial meshes.
\end{example}

\begin{figure}[htb!]
\centering
\subfigure[Time-step ratio $r_n=0.8$]{
\includegraphics[width=2.0in]{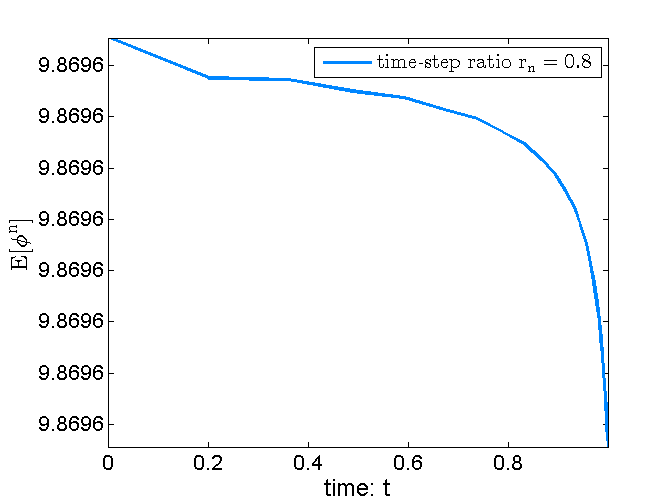}}
\subfigure[Time-step ratio $r_n=1.0$]{
\includegraphics[width=2.0in]{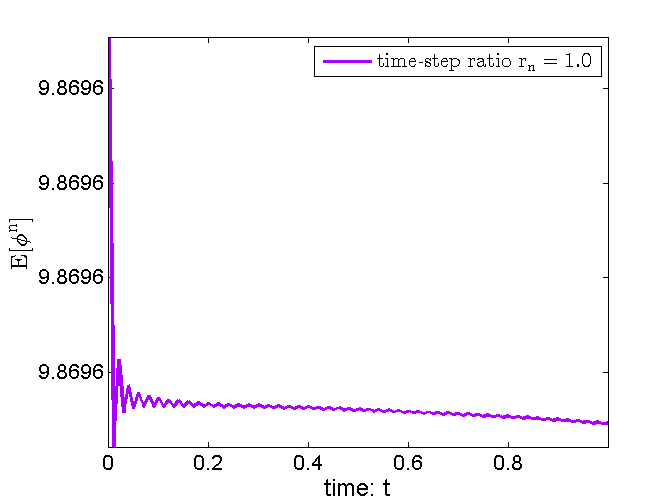}}
\subfigure[Time-step ratio $r_n=1.2$]{
\includegraphics[width=2.0in]{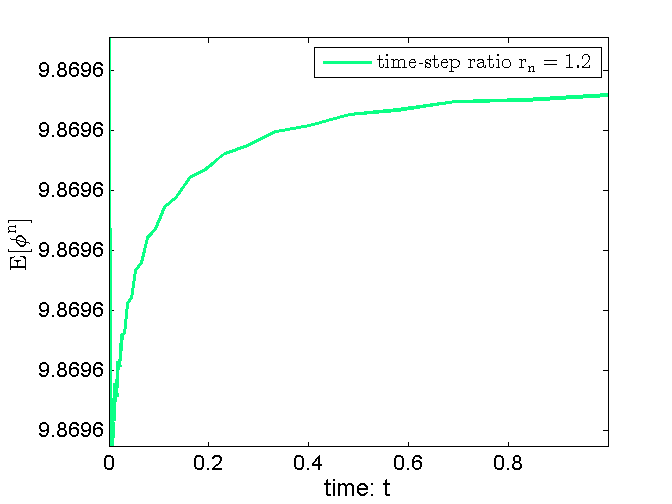}}
\caption{The original energy of the L1$_\mathrm{a}$ scheme
\eqref{scheme:fully implicit L1a TFCH} for  different step ratios $r_n$.}
\label{example:comparison L1a energy with rn}
\end{figure}

At first, we test the discrete original energy $E\kbra{\phi^n}$ defined in
\eqref{def:discrete original energy} using the random initial data,
although no discrete energy law for the L1$_{\mathrm{a}}$
scheme \eqref{scheme:fully implicit L1a TFCH} is built theoretically.
Figure \ref{example:comparison L1a energy with rn} depicts the curves of original energy
for the fractional order $\alpha=0.9$
on the time meshes generated by a fixed
time-step ratio $r_n$ with $N=100$ until time $T=1$.
As observed, the energy dissipation property is violated when the
time-step ratios $r_n\ge 1$, so that the L1$_\mathrm{a}$ scheme
\eqref{scheme:fully implicit L1a TFCH} may be not suitable for practical simulation of the TFCH model.
We thus focus on the numerical computations of the variable-step L1 scheme
\eqref{scheme:fully implicit L1 TFCH} and L1$_\mathrm{h}$ scheme
\eqref{scheme:fully implicit L1h TFCH} in what follows.

\begin{table}[htb!]
\begin{center}
\caption{CPU time  and total time levels with
different adaptive strategies.
}\label{example:compar adap param CPU} \vspace*{0.3pt}
\def\temptablewidth{\textwidth}
{\rule{\temptablewidth}{0.8pt}}
\begin{tabular*}{\temptablewidth}{@{\extracolsep{\fill}}l|cccc}
  \multirow{2}{*}{Time-stepping strategies}
  &\multicolumn{2}{c}{L1 scheme} &\multicolumn{2}{c}{L1$_{\mathrm{h}}$ scheme} \\
                         \cline{2-3}              \cline{4-5}
                         &Total levels &CPU (seconds)
                          &Total levels &CPU (seconds)\\
  \midrule
  uniform step $\tau=5\times10^{-3}$  &6030 &240.526         &6030 &216.339\\
   adaptive steps with $\eta=10$             &487  &28.007          &487  &19.403 \\
   adaptive steps with  $\eta=10^2$           &1092 &51.448          &1089 &39.918  \\
   adaptive steps with  $\eta=10^3$           &3178 &133.169         &3166 &109.450 \\
\end{tabular*}
{\rule{\temptablewidth}{0.5pt}}
\end{center}
\end{table}	

We adopt the graded time meshes $t_k=T_0(k/N_0)^{\gamma}$ together with
the settings $\gamma=3,N_0=30$ and $T_0=0.01$
to resolve the weakly singularity for the THCH model.
The treatment of remainder time interval is a great deal of flexibility
such as the time-stepping strategy below
\cite{Qiao2011An,ZhangQiao:2012AdaptiveCH,HuangYangWei:2020Parallel},
\begin{align}\label{algorithm:adaptive time stepping}
\tau_{ada}
=\max\Bigg\{\tau_{\min},
\frac{\tau_{\max}}{\sqrt{1+\eta\mynormb{\partial_\tau\phi^n}^2}}\Bigg\},
\end{align}
where $\eta>0$ is a user parameter, $\tau_{\max}=0.1$ and $\tau_{\min}=10^{-3}$ are
the predetermined maximum and minimum time steps, respectively.

\begin{figure}[htb!]
\centering
\subfigure[Original energy $E\kbra{\phi^n}$]{
\includegraphics[width=2.0in]{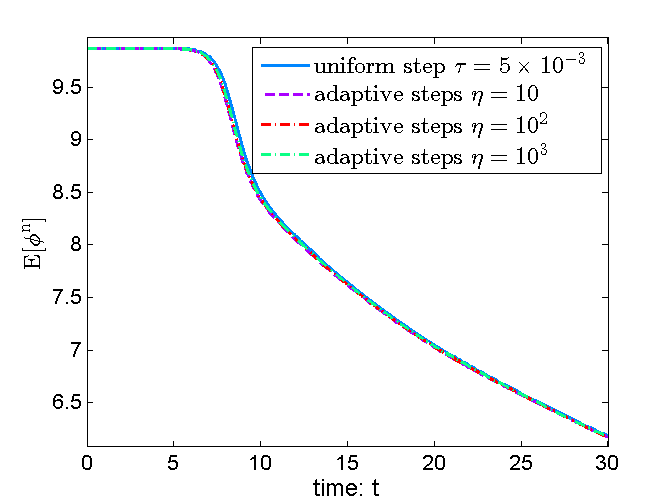}}
\subfigure[Modified energy $\mathcal{E}_\alpha\kbra{\phi^n}$]{
\includegraphics[width=2.0in]{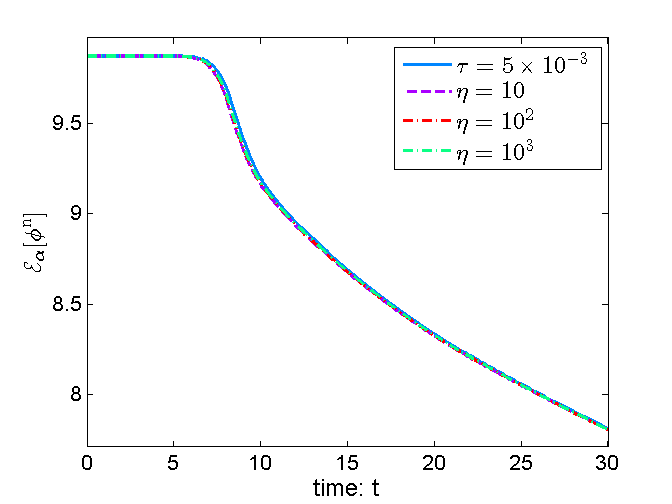}}
\subfigure[Time steps $\tau_n$]{
\includegraphics[width=2.0in]{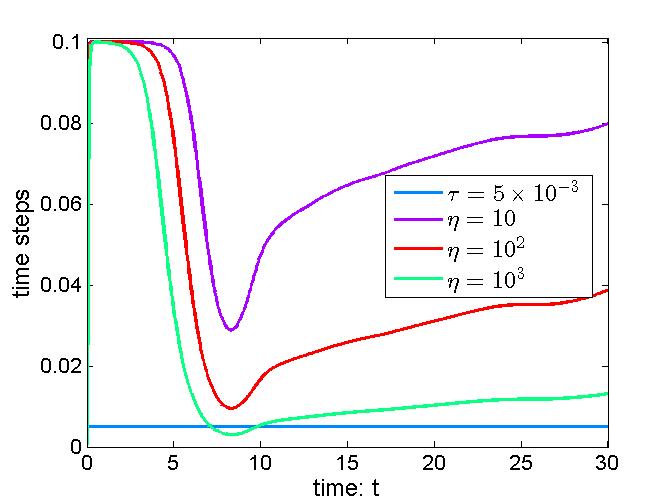}}
\caption{Energy curves by uniform step and adaptive strategy with different parameters $\eta$.}
\label{example:compar adap param}
\end{figure}

\begin{figure}[htb!]
\centering
\subfigure[The profile of $\phi$ with fractional order $\alpha=0.4$ at time $t=30,100,300,500.$]{
\includegraphics[width=1.47in]{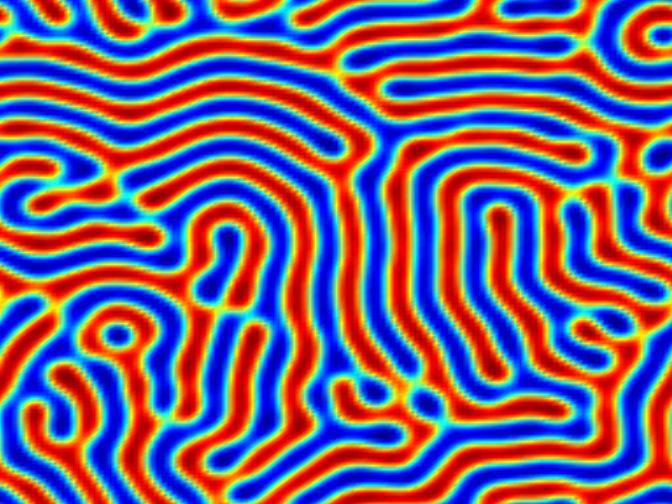}
\includegraphics[width=1.47in]{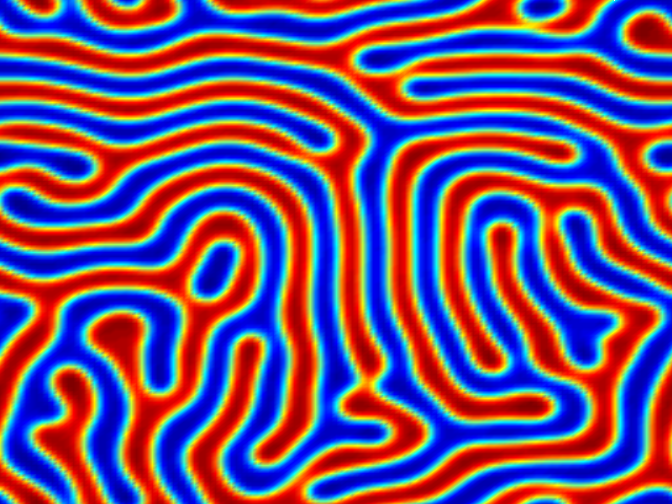}
\includegraphics[width=1.47in]{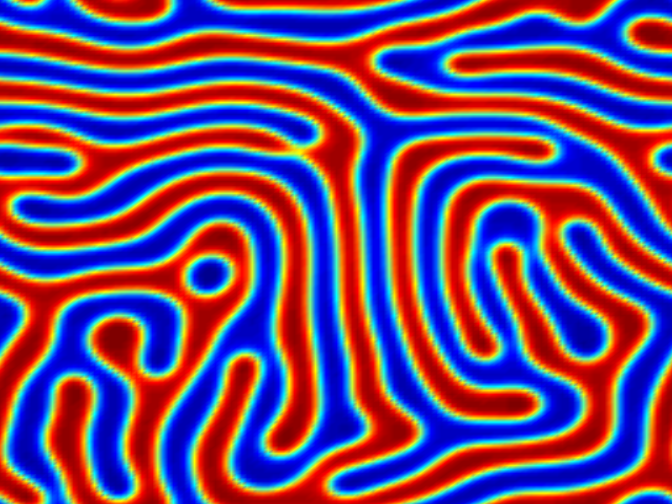}
\includegraphics[width=1.47in]{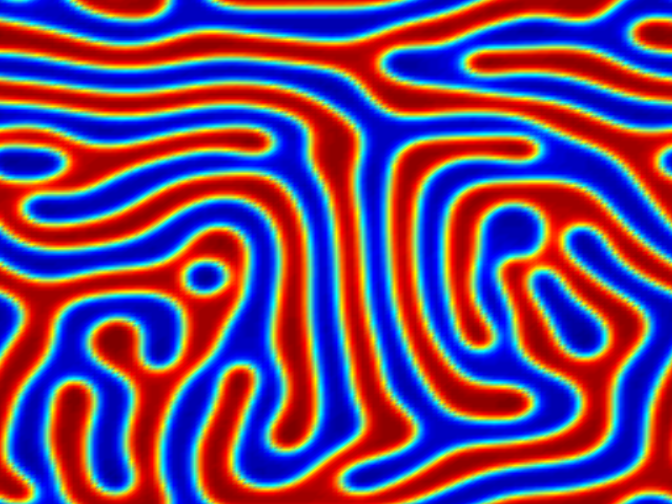}}\\
\subfigure[The profile of $\phi$ with fractional order $\alpha=0.7$ at time $t=30,100,300,500.$]{
\includegraphics[width=1.47in]{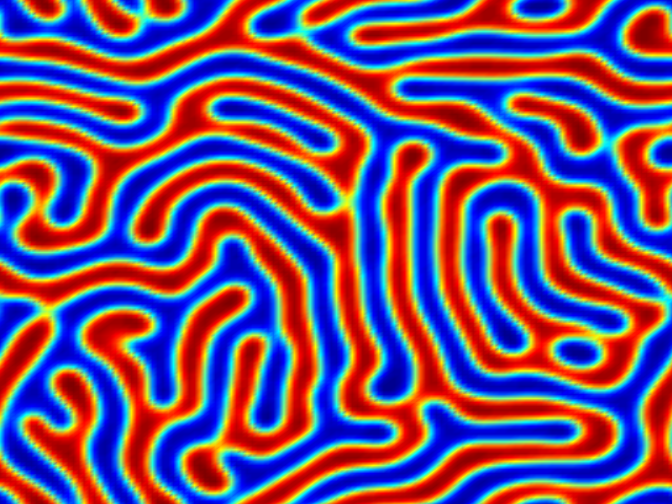}
\includegraphics[width=1.47in]{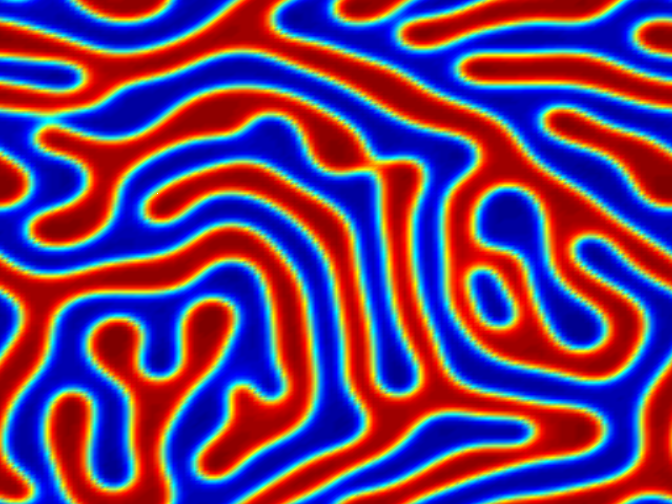}
\includegraphics[width=1.47in]{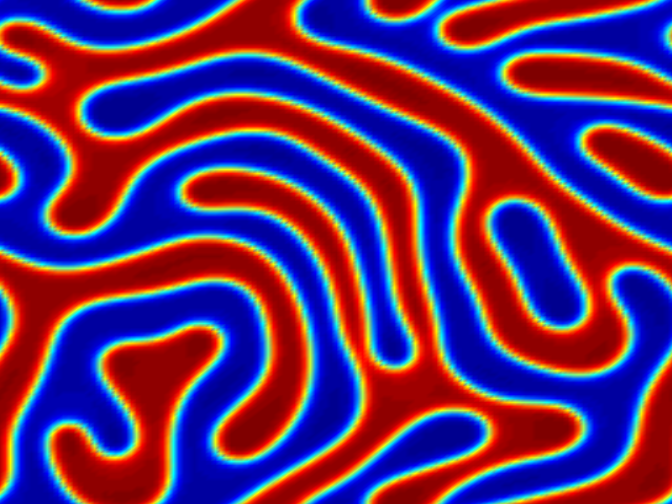}
\includegraphics[width=1.47in]{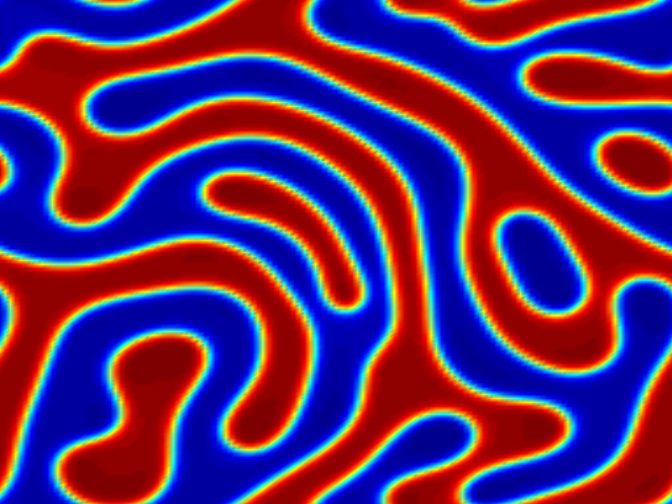}}\\
\subfigure[The profile of $\phi$ with fractional order $\alpha=0.9$ at time $t=30,100,300,500.$]{
\includegraphics[width=1.47in]{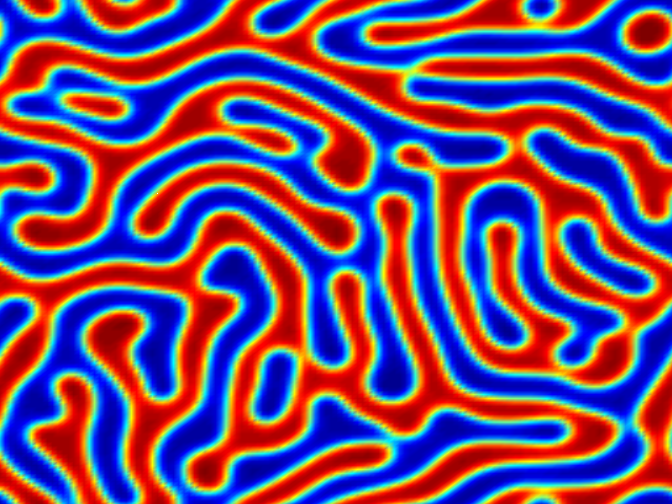}
\includegraphics[width=1.47in]{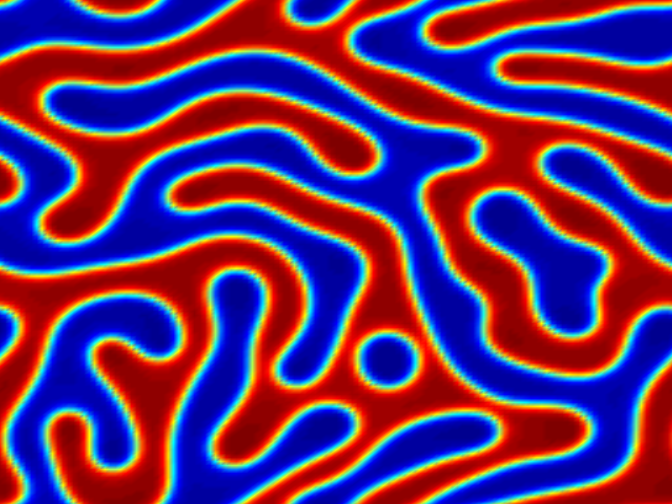}
\includegraphics[width=1.47in]{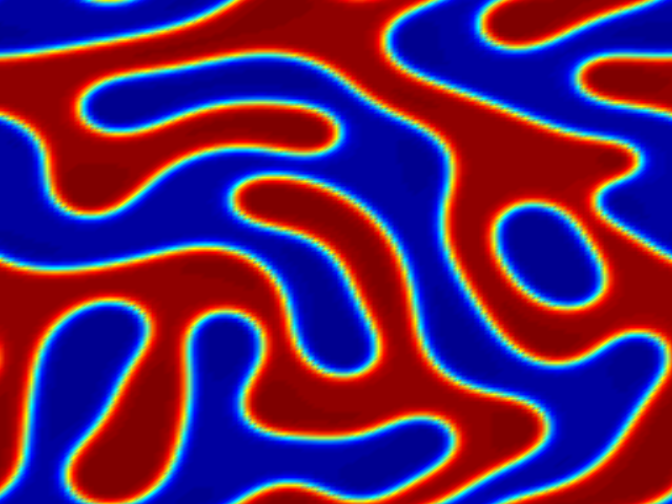}
\includegraphics[width=1.47in]{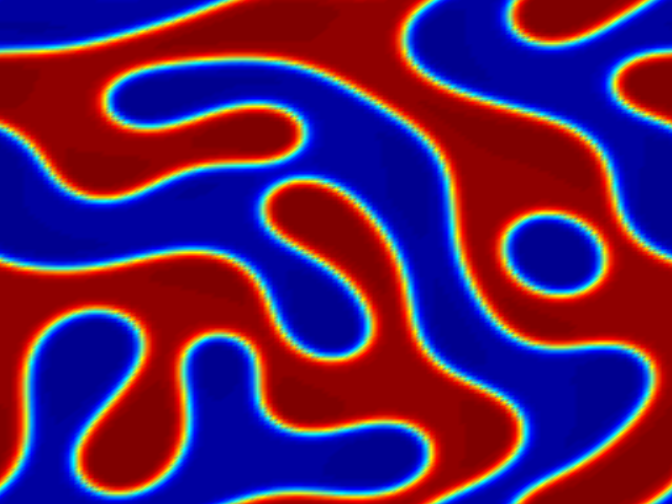}}\\
\caption{Snapshots of dynamic coarsening processes for different fractional orders $\alpha$.}
\label{example:FCH coarsen snap}
\end{figure}

To test the numerical performance of the adaptive time-stepping algorithm
\eqref{algorithm:adaptive time stepping},
we perform a comparative study by running
the L1 scheme \eqref{scheme:fully implicit L1 TFCH}
and the L1$_{\mathrm{h}}$ scheme \eqref{scheme:fully implicit L1h TFCH}
on different time steps.
We first apply a small uniform time step $\tau=5\times10^{-3}$
to obtain the reference solution.
Then we repeat the numerical simulation by
using the adaptive time-stepping strategy  with three different parameters
$\eta=10,10^2,10^3$, respectively.

\begin{figure}[htb!]
\centering
\subfigure[Original energy $E\kbra{\phi^n}$]{
\includegraphics[width=2.0in]{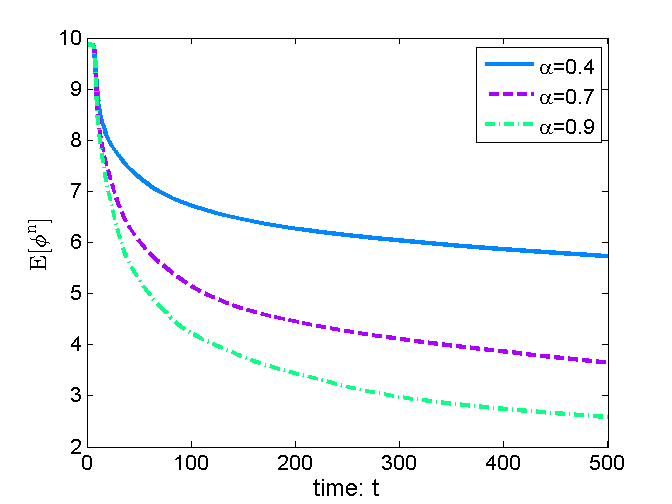}}
\subfigure[Variational energy $\mathcal{E}_\alpha\kbra{\phi^n}$]{
\includegraphics[width=2.0in]{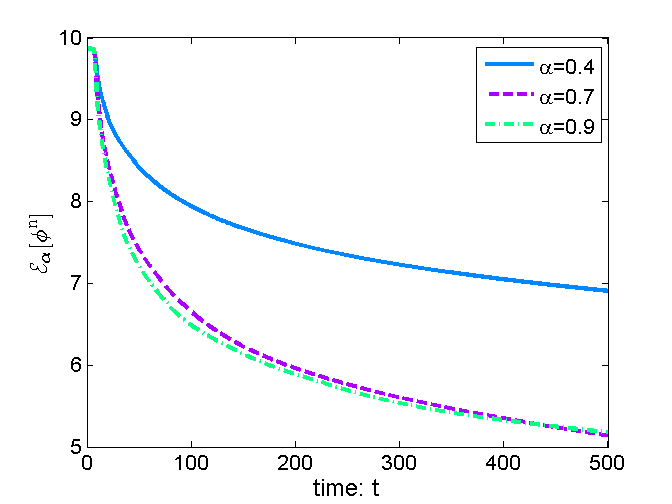}}
\subfigure[Adaptive time steps $\tau_n$]{
\includegraphics[width=2.0in]{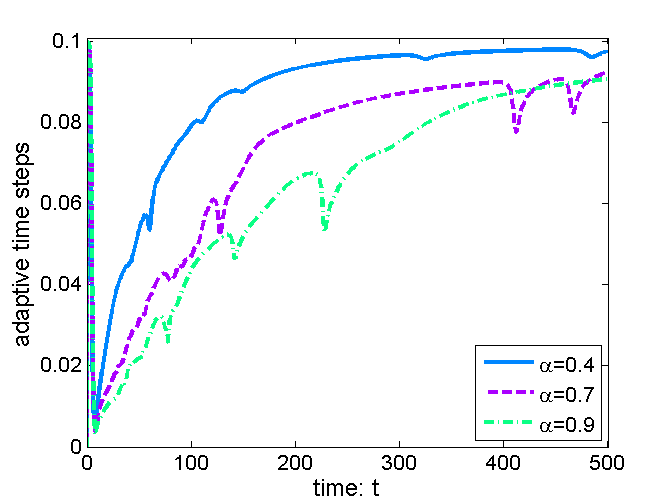}}
\caption{Numerical results of  the TFCH model with different fractional orders $\alpha$.}
\label{example:FCH coarsen energy}
\end{figure}

The numerical results are summarized in Figure \ref{example:compar adap param}.
As can be seen, the numerical results
using adaptive time-stepping are comparable to the reference solution.
Also, one can observe that the adaptive time-steps are adjusted promptly by
the parameters $\eta$: large (small) $\eta$
reinforces (reduces) the restriction to the time step sizes.
The corresponding CPU time (in seconds) and the total time levels
for different time-stepping strategies are listed in Table \ref{example:compar adap param CPU}.
The effectiveness of the adaptive time-stepping algorithm
makes the long-time dynamics simulation practical.
Note that, the numerical results of the L1$_\mathrm{h}$ scheme \eqref{scheme:fully implicit L1h TFCH} are
quite similar to those of the L1 scheme \eqref{scheme:fully implicit L1 TFCH}, and we thus omit them for brevity.

Finally, we perform the numerical simulation by using
the adaptive time-stepping strategy \eqref{algorithm:adaptive time stepping}
with the parameter $\eta=10^3$ until time $T=500$.
The rest settings are kept the same as in the previous example.
The profile of $\phi$ for the TFCH model \eqref{cont:TFCH model}
with different fractional orders $\alpha$ are
depicted in Figure \ref{example:FCH coarsen snap}.
They are consistent with the coarsening dynamics process reported in \cite{LiuChengWangZhao:2018TimeAC-CH,ZhaoChenWang:2019OnPowerLaw}.
The evolutions of the numerical energies and adaptive time steps
during the coarsening dynamics are
depicted in Figure \ref{example:FCH coarsen energy}.
They suggest that the proposed variable-step methods effectively
capture the multiple time scales in the long-time dynamical simulations.

\lan{
\section*{Acknowledgements}
The authors would like to thank the editor and the anonymous referees for their valuable suggestions and some recent works on the positive definiteness of quadratic form with discrete convolution kernels. They are helpful in improving the quality of the paper.  
}


\end{document}